\documentclass[reqno,fleqn,onecolumn]{amsproc}%
\usepackage{amsfonts}
\usepackage{amsmath}
\usepackage{amssymb}
\usepackage{graphicx}
\usepackage[none]{hyphenat}%
\theoremstyle{plain}

\numberwithin{equation}{section}
\begin{document}
\title{}

\begin{center}
{\LARGE \textbf{Game pricing and double sequence of random variables}}

\bigskip

\bigskip

\textbf{Yukio Hirashita}

\bigskip

\textit{Faculty of Liberal Arts, Chukyo University, Nagoya, Aichi 466-8666,
Japan}

\bigskip
\end{center}

\noindent-----------------------------------------------------------------------------------------------------------

\noindent{\large \textbf{Abstract}}

In this paper, we study a game with positive or plus infinite expectation and
determine the optimal proportion of investment for maximizing the limit
expectation of growth rate per attempt. With this objective, we introduce a
new pricing method in which the price is different from that obtained by the
Black-Scholes formula for a European option.



\noindent-----------------------------------------------------------------------------------------------------------

\bigskip

\bigskip

\noindent{\large \textbf{1. Introduction}}

The portfolio pricing equation (Luenberger (1998) 9.7) is useful for
determining prices only if the optimal portfolio is already known. In this
paper, we determine both the price and optimal proportion of investment for
any effective game (Section 6).

The determination of the utility function is more experimental than
mathematical. In general, despite the equality $(\Pi_{j=1}^{n}X_{j}%
)^{1/n}=\exp(\Sigma_{j=1}^{n}\log X_{j}/n)$, the two expectations
$E[(\Pi_{j=1}^{n}X_{j})^{1/n}]$ and $\exp(E[\Sigma_{j=1}^{n}\log X_{j}/n])$
are not equal for a sequence of independent random variables $\{X_{j}>0$ $|$
$j$ $=1,2,\cdots,n\}$. Therefore, we use neither the notion of utility from
consumption nor the law of large numbers (Luenberger (1998) 15.2).

The investor should repeatedly invest a fixed proportion of his or her own
current capital without borrowing. As a rule, if the investor invests $1$
dollar, then he or she receives $a(x)$ dollars (including the invested $1$
dollar) with a distribution function $F(x)$ defined on an interval $I$. For
simplicity, we omit the currency notation. Let $M>0$ be the investor's
capital, $u>0$ the price of the game, and $0\leq t\leq1$ the proportion of
investment. Then, after one attempt, he or she has capital of
$Mta(x)/u+M(1-t)$ if $x$ occurs. It should be noted that the reserved part
$M(1-t)$ does not include the interest, that is the custom, for example, in
foreign exchange accounts.

Let $M_{n}>0$ be the capital after $n$ attempts. In general, growth rate
implies $M_{n+1}/M_{n}-1$ or $\log(M_{n+1}/M_{n})$ after one attempt. However,
for the purposes of succinctness in this paper $M_{n+1}/M_{n}$ is used to
define the growth rate. In this context, the growth rate per attempt is
defined as $\left(  M_{n}/M_{0}\right)  ^{\frac{1}{n}}$.

Without dealing with $\left(  M_{n}/M_{0}\right)  ^{\frac{1}{n}}$ directly,
this paper defines a double sequence of random variables $\{X_{N,n}\}$ with
respect to the bounded step functions $\{f_{N}(x)\}$ such that $\lim
_{N\rightarrow+\infty}$ $f_{N}(x)=a(x)$ (Section 5). It is shown that the
finite limit $\lim_{\substack{n\rightarrow+\infty\\N\rightarrow+\infty
}}E[X_{N,n}]$ exists if, and only if, the game is effective. In this case, the
equalities $\lim_{\substack{n\rightarrow+\infty\\N\rightarrow+\infty
}}E[X_{N,n}]=G_{u}(t):=\exp(\int_{I}\log(a(x)t/u-t+1)d(F(x)))$ and
$\lim_{\substack{n\rightarrow+\infty\\N\rightarrow+\infty}}V[X_{N,n}]=0$ are
obtained. These equalities again support the well-known assertion that
although in principle an investor may choose any utility function, a
repetitive situation tends to hammer the utility into one that is close to the
logarithm (Luenberger (1998) 15.4, Kelly (1956)).

We study the optimal proportion of investment, $t_{u}$, for the price $u>0$ in
order to maximize the limit expectation of growth rate per attempt. In order
to determine the price of the game, we require a riskless interest rate,
$r>0,$ for a particular period. The equation $G_{u}(t_{u})=r+1$ (if $r$ is
simple) or $G_{u}(t_{u})=e^{r\text{ }}$(if $r$ is continuously compounded)\ is
used to determine the price of a game. If $a(x)\geq0$ for each $x\in I,$ then
the existence and uniqueness of the price are guaranteed by the fact that
$G_{u}(t_{u})$ is continuous and strictly decreases from $+\infty$ to $1$ with
respect to $0<u<E:=\int_{I}a(x)d(F(x))$ (Theorem 4.1). In this context, the
price of the St. Petersburg game (Daniel Bernoulli (1738; English trans.
1954)) is determined to be $5.1052$ if the riskless interest rate is $4\%$
(Example 6.4). On the other hand, the Black-Scholes formula is deduced from
the equation $E/u=e^{r},$ where $E$ is the expectation of a European option
(Example 6.6).

\bigskip

\bigskip

\noindent{\large \textbf{2. Optimal proportion of investment }}

\bigskip

Assume that the profit function $a(x)$ is measurable with the distribution
function $F(x)$ defined on an interval $I\subseteq(-\infty$, $+\infty)$. Set
$\xi:=\inf_{x\in I}\,a(x)$. We also assume that $\xi>-\infty$ and $\xi$ is the
essential infimum of $a(x)$, that is, $\int_{a(x)<\xi+\varepsilon}d(F(x))$
$>0$ for each $\varepsilon>0$. Further, assume that $a(x)$ is not a constant
function (a.e.), that is, $\int_{a(x)<\xi+\delta}d(F(x))<1$ for some
$\delta>0.$

We use the following notation.
\begin{equation}
E:=\int_{I}a(x)d(F(x)),\quad H:=\int_{I}\frac{1}{a(x)}d(F(x)),\quad H_{\xi
}:=\int_{I}\frac{1}{a(x)-\xi}d(F(x)). \tag{1}%
\end{equation}
In this paper, we assume that $E>0$. If $\int_{a(x)=\xi}d(F(x))>0$, we define
$H_{\xi}$ $=+\infty$ and $1/H_{\xi}$ $=0.$

Since $a(x)$ is not constant, we have $\xi<E$, $H_{\xi}>0$, and $1/H_{\xi
}<+\infty$. From the relation%
\begin{align*}
1  &  =\left(  \int_{I}\sqrt{a(x)-\xi}\times\frac{1}{\sqrt{a(x)-\xi}%
}d(F(x))\right)  ^{2}\\
&  <\int_{I}(a(x)-\xi)d(F(x))\times\int_{I}\frac{1}{a(x)-\xi}d(F(x))=(E-\xi
)H_{\xi},
\end{align*}
we have $\xi+1/H_{\xi}<E$. In particular, if $\xi=0,$ $0\leq1/H<E.$ If $\xi
>0$, then using $1/\xi\geq1/a(x)$ and $1=\sqrt{a(x)}\times(1/\sqrt{a(x)})$, we
have $\xi<1/H<E.$

\bigskip

For price $u>0$, let $t_{u}\in\lbrack0,$ $1]$ be the optimal proportion of
investment. The precise definition of the term "optimal" and its significance
is provided in Section 5. Here, we present certain properties of $t_{u}$ in
order to explain the approximate outline of the paper.

\bigskip

\noindent(\textbf{a}) If $u>E,$ $t_{u}=0$.

Assume that $u>E$ and $t\in(0$, $1]$, then the expectation of profits,
$Mt\int_{I}a(x)/u$ $d(F(x))+M(1-t)=M-M(1-E/u)t$, is less than $M$. More
precisely, using Jensen's inequality, we have $G_{u}(t)<1-(1-E/u)t<1=G_{u}(0)$
for each $t\in(0$, $1]$. Therefore, $t_{u}=0$.

In the proof of Theorem 5.1, we will show that:
\begin{equation}
\{u\text{ }|\text{ }t_{u}=0\}=\left\{
\begin{array}
[c]{cc}%
\lbrack E,+\infty), & \text{if }E<+\infty,\\
\phi, & \text{if }E=+\infty\text{.}%
\end{array}
\right.  \tag{2}%
\end{equation}

\bigskip

\noindent(\textbf{b}) If $\xi>0$ and $0<u\leq\xi,$ then $t_{u}=1$.

From $0<u\leq\xi\leq a(x)$ and $t\in\lbrack0$, $1)$, we have $Mta(x)/u+M(1-t)$
$=Ma(x)/u-M(1-t)(a(x)/u-1)\leq Ma(x)/u$ for each $x\in I.$ This implies that
$G_{u}(t)<G_{u}(1)$ for each $t\in\lbrack0$, $1),$ that is, $t_{u}=1$.

Accordingly, in the proof of Theorem 5.1, we will also show that
\begin{equation}
\{u\text{ }|\text{ }t_{u}=1\}=\left\{
\begin{array}
[c]{cc}%
(0\text{, }1/H], & \text{if }\xi>0,\text{ or }\xi=0\text{ and }H<+\infty,\\
\phi, & \text{if }\xi<0,\text{ or }\xi=0\text{ and }H=+\infty,
\end{array}
\right.  \tag{3}%
\end{equation}
which yields a maximum price of $1/H$ at which all the capital should be
repeatedly invested.

\bigskip

\noindent(\textbf{c}) If $\max(0,\xi)<u$, $t_{u}\leq u/(u-\xi)$.

If $t>u/(u-\xi),$ $u-\xi-u/t>0$. Therefore, the negative result $Mta(x)/u$
$+M(1-t)<0$ occurs with a positive probability $\int_{\xi\leq a(x)<\xi
+(u-\xi-u/t)}d(F(x))>0.$ This contradicts the concept of continual investment
without borrowing.

In the proof of Theorem 5.1, the existence of $t_{u}$ is shown such that:
\begin{equation}
\{u\text{ }|\text{ }t_{u}=\frac{u}{u-\xi}\}=\left\{
\begin{array}
[c]{cc}%
(0\text{, }\xi+1/H_{\xi}], & \text{if }\xi\leq0\text{ and }\xi+1/H_{\xi}>0,\\
\phi, & \text{if }\xi>0\text{ or }\xi+1/H_{\xi}\leq0.
\end{array}
\right.  \tag{4}%
\end{equation}

\bigskip

\noindent\textbf{( d )} Theorem 5.1 also shows that $t_{u}\notin\{0,$ $1,$
$u/(u-\xi)\}$ if and only if $1/H$ $<u<E$ (if $\xi\geq0$) or $\max(0$,
$\xi+1/H_{\xi})<u<E$ (if $\xi<0$). In this case, $t_{u}$ can be uniquely
determined by the property:%
\begin{equation}
\int_{I}\frac{a(x)-u}{a(x)t_{u}-ut_{u}+u}d(F(x))=0. \tag{5}%
\end{equation}

\bigskip

\bigskip

\noindent{\large \textbf{3. Pre-optimal proportion}}

\bigskip

We denote the integral $\int_{I}(a(x)-\beta)/(a(x)z-z\beta+\beta)d(F(x))$ by
$w_{\beta}(z)$, in which $z$ and $\beta$ are complex variables.

\bigskip

\noindent\textbf{Lemma 3.1.} \textit{The function }$w_{\beta}(z)$\textit{ is
holomorphic with respect to two complex variables }$z:=t+si$\textit{ and
}$\beta:=u+hi$\textit{ such that,}

(a)\textit{ }$\max(\varepsilon$\textit{, }$\xi)<u<L,$

(b)\textit{ }$|h|<\varepsilon^{6}/(32(L+1)R^{2}),$

(c)\textit{ }$|z|<R$\textit{ and }$z\notin\{\left\vert s\right\vert
\leq\varepsilon\}\cap\{t\leq\varepsilon$\textit{ \ or }$t\geq u/(u-\xi
)-\varepsilon\},$

\textit{where }$0<\varepsilon<\min(1/2$\textit{, }$u/(2(u-\xi)))$\textit{,
}$\max(\varepsilon$\textit{, }$\xi)<L<+\infty$\textit{, }$\max(2$, $u/(u$
$-\xi))<R<+\infty$\textit{, }$i:=\sqrt{-1}$, $\operatorname{Im}(z)=s$
\textit{and} $\operatorname{Im}(\beta)=h$.

\bigskip

\noindent\textbf{Proof. }We obtain certain operator exchange properties such
as
\[
\frac{\partial}{\partial t}w_{\beta}(z)=\int_{I}\frac{\partial}{\partial
t}\left(  \frac{a(x)-\beta}{a(x)z-z\beta+\beta}\right)  d(F(x))
\]
by proving that the related integrands are bounded. Because $\left(
a(x)-\beta\right)  /(a(x)z$ $-z\beta+\beta)$ satisfies the Cauchy-Riemann
equations, $w_{\beta}(z)$ is shown to be holomorphic due to Hartogs's theorem.

It should be noted that the condition (a) above leads to $\beta\neq0$, and if
$a(x)\neq\beta,$ then we have%
\[
\frac{a(x)-\beta}{a(x)z-z\beta+\beta}=\frac{1}{z-\frac{1}{1-\frac{a(x)}{\beta
}}}.
\]
In the following four cases, we assume that $a(x)\neq\beta.$

In this proof, we will frequently use the inequality $\left\vert
1/(1-z)\right\vert \leq\left\vert 2/z\right\vert $ if $|z|\geq2.$%

$<$%
Case 1%
$>$%
\textbf{ }$|a(x)|\geq8(L+1)/\varepsilon.$

As a result of the conditions, we have $\left\vert a(x)/u\right\vert
>\left\vert a(x)/L\right\vert >\left\vert a(x)/(L+1)\right\vert $
$\geq8/\varepsilon$ $>16,$ which leads to $\left\vert 1/(1-a(x)/u)\right\vert
\leq\left\vert 2/(a(x)/u))\right\vert <\left\vert 2(L+1)/a(x)\right\vert $
$\leq\varepsilon/4$. On the other hand, the inequality $\left\vert
a(x)/\beta\right\vert >\left\vert \frac{a(x)}{L+1}\right\vert \geq
8/\varepsilon>16$ leads to $\left\vert 1/(1-a(x)/\beta)\right\vert
\leq\left\vert 2/(a(x)/\beta)\right\vert <\left\vert 2(L+1)/a(x)\right\vert
\leq\varepsilon/4$, where $\left\vert \beta\right\vert \leq u+\left\vert
h\right\vert $ $<L+1.$

Moreover, from $\xi\leq a(x)$ we have $1-a(x)/u\leq(u-\xi)/u.$ If $1-a(x)/u>0$
then $u/(u-\xi)\leq1/(1-a(x)/u)$, which leads to $\left\vert
z-1/(1-a(x)/u))\right\vert >\varepsilon$ due to (c). If $1-a(x)/u<0$ then
$1/(1-a(x)/u)<0$, which leads to $\left\vert z-1/(1-a(x)/u))\right\vert
>\varepsilon$ due to (c). If $1-a(x)/u=0$ then $L>u=|a(x)|\geq
8(L+1)/\varepsilon$, which is a contradiction. Therefore, we have%
\[
\left\vert z-\frac{1}{1-\frac{a(x)}{\beta}}\right\vert =\left\vert z-\frac
{1}{1-\frac{a(x)}{u}}+\frac{1}{1-\frac{a(x)}{u}}-\frac{1}{1-\frac{a(x)}{\beta
}}\right\vert >\frac{\varepsilon}{2},
\]
which establishes%
\[
\frac{1}{\left\vert z-\frac{1}{1-\frac{a(x)}{\beta}}\right\vert }<\frac
{2}{\varepsilon}.
\]
Moreover, using $\left\vert a(x)/(1-a(x)/\beta)\right\vert \leq\left\vert
2\beta\right\vert $, $\left\vert 1/(1-a(x)/\beta)\right\vert <\varepsilon/4$,
and $1/\left\vert \beta\right\vert $ $<1/\varepsilon$, we have%
\[
\left\vert \frac{a(x)}{\beta^{2}\left(  z-\frac{1}{1-\frac{a(x)}{\beta}%
}\right)  ^{2}\left(  1-\frac{a(x)}{\beta}\right)  ^{2}}\right\vert <\frac
{2}{\varepsilon^{2}}.
\]

\bigskip%

$<$%
Case 2%
$>$%
\textbf{ }$|a(x)|<8(L+1)/\varepsilon$ and $\left\vert a(x)/\beta-1\right\vert
\leq\varepsilon/R.$

Since $\left\vert \left(  a(x)/\beta-1\right)  z\right\vert \leq\varepsilon$,
we have $\left\vert a(x)z/\beta-z+1\right\vert \geq1-\varepsilon$. Therefore,
\[
\left\vert \frac{\frac{a(x)}{\beta}-1}{\frac{a(x)}{\beta}z-z+1}\right\vert
\leq\frac{\frac{\varepsilon}{R}}{1-\varepsilon}<\frac{2\varepsilon}{R}.
\]
Moreover, using $1/\left\vert \beta\right\vert <1/\varepsilon$ and
$1/(1-\varepsilon)<2,$ we have%
\[
\left\vert \frac{a(x)}{\beta^{2}\left(  z-\frac{1}{1-\frac{a(x)}{\beta}%
}\right)  ^{2}\left(  1-\frac{a(x)}{\beta}\right)  ^{2}}\right\vert
=\left\vert \frac{a(x)}{\beta^{2}\left(  \frac{a(x)}{\beta}z-z+1\right)  ^{2}%
}\right\vert <\frac{32(L+1)}{\varepsilon^{3}}.
\]

\bigskip%

$<$%
Case 3%
$>$
$|a(x)|<8(L+1)/\varepsilon$, $\left\vert a(x)/\beta-1\right\vert
>\varepsilon/R,$ and $\left\vert a(x)/u-1\right\vert $ $>\varepsilon/(2R).$

From $1/\left\vert \beta\right\vert <1/\varepsilon,$ $1/u<1/\varepsilon$, and
condition (b) mentioned above, we have
\[
\left\vert \frac{1}{1-\frac{a(x)}{u}}-\frac{1}{1-\frac{a(x)}{\beta}%
}\right\vert =\left\vert \frac{a(x)hi}{u\beta(1-\frac{a(x)}{u})(1-\frac
{a(x)}{\beta})}\right\vert <\frac{\varepsilon}{2}.
\]
Therefore, as in Case 1, we obtain%
\[
\frac{1}{\left\vert z-\frac{1}{1-\frac{a(x)}{\beta}}\right\vert }<\frac
{2}{\varepsilon}.
\]
Moreover, this implies that%
\[
\left\vert \frac{a(x)}{\beta^{2}\left(  z-\frac{1}{1-\frac{a(x)}{\beta}%
}\right)  ^{2}\left(  1-\frac{a(x)}{\beta}\right)  ^{2}}\right\vert
<\frac{32(L+1)R^{2}}{\varepsilon^{7}}.
\]

\bigskip%

$<$%
Case 4%
$>$
$|a(x)|<8(L+1)/\varepsilon$, $\left\vert a(x)/\beta-1\right\vert
>\varepsilon/R,$ and $\left\vert a(x)/u-1\right\vert \leq\varepsilon/(2R).$

This case is void as shown below:%
\[
\frac{\varepsilon}{2R}<\left\vert \frac{a(x)}{\beta}-1\right\vert -\left\vert
\frac{a(x)}{u}-1\right\vert \leq\left\vert \frac{a(x)}{\beta}-\frac{a(x)}%
{u}\right\vert =\frac{\left\vert ha(x)\right\vert }{\left\vert \beta
u\right\vert }<\frac{\varepsilon^{3}}{4R^{2}},
\]
which leads to the contradiction $4<2R<\varepsilon^{2}<1/4.$

\bigskip%

$<$%
Conclusion%
$>$
From the inequalities mentioned above, the four integrands on the right-hand
side of the following equalities are bounded. Therefore, the Cauchy-Riemann
equations for $w_{\beta}(z)$ hold.%
\begin{align}
\frac{\partial}{\partial t}w_{\beta}(z)  &  =\int_{a(x)\neq\beta}\frac
{-1}{\left(  z-\frac{1}{1-\frac{a(x)}{\beta}}\right)  ^{2}}d(F(x)),\tag{6}\\
\frac{\partial}{\partial s}w_{\beta}(z)  &  =\int_{a(x)\neq\beta}\frac
{-i}{\left(  z-\frac{1}{1-\frac{a(x)}{\beta}}\right)  ^{2}}d(F(x)),\nonumber\\
\frac{\partial}{\partial u}w_{\beta}(z)  &  =\int_{a(x)\neq\beta}\frac
{-a(x)}{\beta^{2}\left(  \frac{a(x)}{\beta}z-z+1\right)  ^{2}}%
d(F(x)),\nonumber\\
\frac{\partial}{\partial h}w_{\beta}(z)  &  =\int_{a(x)\neq\beta}\frac
{-a(x)i}{\beta^{2}\left(  \frac{a(x)}{\beta}z-z+1\right)  ^{2}}%
d(F(x)).\nonumber
\end{align}

\hfill$\square$

\bigskip

Henceforth, in this Section, we assume that $\max(0$, $\xi)<u<E$ and $0<t$
$<u/(u-\xi)$. It should be noted that $a(x)t-ut+u\geq\xi t-ut+u=(u-\xi
)(u/(u-\xi)-t)>0$ for each $x\in I.$

\bigskip

\noindent\textbf{Lemma 3.2}. $w_{u}(t)$\textit{ is strictly decreasing with
respect to }$t$.

\bigskip

\noindent\textbf{Proof. }According to Lemma 3.1, we have
\[
\frac{\partial}{\partial t}w_{u}(t)=-\int_{I}\left(  \frac{a(x)-u}%
{a(x)t-ut+u}\right)  ^{2}d(F(x))<0.
\]

\hfill$\square$

\bigskip

\noindent\textbf{Lemma 3.3}. $\lim_{t\rightarrow0^{+}}w_{u}(t)=E/u-1$.

\bigskip

\noindent\textbf{Proof. }Since $(a(x)-u)/(a(x)t-ut+u)$ is strictly decreasing
($a(x)\neq u$) with respect to $0<t<u/(u-\xi)$, using Lebesgue (monotone
convergence) theorem, we obtain%
\[
\lim_{t\rightarrow0^{+}}w_{u}(t)=\lim_{t\rightarrow0^{+}}\int_{I}\frac
{a(x)-u}{a(x)t-ut+u}d(F(x))=\int_{I}\frac{a(x)-u}{u}d(F(x))=\frac{E}{u}-1.
\]

\hfill$\square$

\bigskip

\noindent\textbf{Lemma 3.4. }$\lim_{t\rightarrow(u/(u-\xi))^{-}}%
w_{u}(t)=(1-\xi/u)H_{\xi}(\xi+1/H_{\xi}-u).$

\bigskip

\noindent\textbf{Proof. }Using the same principle as above, we obtain%
\begin{align*}
\lim_{t\rightarrow(u/(u-\xi))^{-}}w_{u}(t)  &  =\frac{u-\xi}{u}\int_{I}%
\frac{a(x)-u}{a(x)-\xi}d(F(x))\\
&  =\frac{u-\xi}{u}(1-(u-\xi)H_{\xi})=(1-\xi/u)H_{\xi}(\xi+1/H_{\xi}-u).
\end{align*}

\hfill$\square$

\bigskip

From the above lemmas, if $\max(0,\xi+1/H_{\xi})<u<E$, then $\lim
_{t\rightarrow0^{+}}w_{u}(t)>0$ and $\lim_{t\rightarrow(u/(u-\xi))^{-}}%
w_{u}(t)<0$. Thus, the equation $w_{u}(t)=0$ has the only solution
$\widetilde{t}_{u}\in(0$, $u/(u-\xi))$, and we refer to it as
\textit{pre-optimal proportion}. Note that, due to Lemma 3.1 and the inverse
mapping theorem, $\widetilde{t}_{u}$ is continuous with respect to $u$.

\bigskip

\noindent\textbf{Lemma 3.5}. \textit{If }$\xi>0,$\textit{ }$\xi+1/H_{\xi
}<1/H.$

\bigskip

\noindent\textbf{Proof. }Since $w_{u}(t)$ is strictly decreasing and
$1<u/(u-\xi)$, $w_{u}(1)>\lim_{t\rightarrow(u/(u-\xi))^{-}}$ $w_{u}(t)$, that
is,%
\[
\int_{I}\frac{a(x)-u}{a(x)-u+u}d(F(x))=H(1/H-u)>(1-\xi/u)H_{\xi}(\xi+1/H_{\xi
}-u)\text{ }%
\]
for each $\xi<u<E.$ If $1/H\leq\xi+1/H_{\xi}$, then selecting $u=\xi+1/H_{\xi
}<E$ leads to the contradiction that $0\geq H(1/H-u)>(1-\xi/u)H_{\xi}%
(\xi+1/H_{\xi}-u)=0$.\hfill$\square$

\bigskip

\noindent\textbf{Lemma 3.6.}\textit{ If }$\xi\geq0,$\textit{ then }%
$\widetilde{t}_{u}$\textit{ is strictly decreasing with respect to }$u\in
(\xi+1/H_{\xi},$ $E)$.

\bigskip

\noindent\textbf{Proof.} From Lemma 3.1, $\widetilde{t}_{u}$ is analytic.
Using
\[
\frac{\partial}{\partial u}\left(  \frac{a(x)-u}{a(x)\widetilde{t}%
_{u}-u\widetilde{t}_{u}+u}\right)  =-\frac{a(x)+(a(x)-u)^{2}\frac
{d\widetilde{t}_{u}}{du}}{(a(x)\widetilde{t}_{u}-u\widetilde{t}_{u}+u)^{2}},
\]
$dw_{u}(\widetilde{t}_{u})/du=0,$ and $a(x)\geq\xi\geq0$, we obtain%
\[
\frac{d\widetilde{t}_{u}}{du}=\frac{-\int_{I}\frac{a(x)}{(a(x)\widetilde
{t}_{u}-u\widetilde{t}_{u}+u)^{2}}d(F(x))}{\int_{I}\frac{(a(x)-u)^{2}%
}{(a(x)\widetilde{t}_{u}-u\widetilde{t}_{u}+u)^{2}}d(F(x))}<0.
\]

\hfill$\square$

\bigskip

\noindent\textbf{Lemma 3.7.} $\lim_{u\rightarrow E^{-}}\widetilde{t}_{u}=0.$

\bigskip

\noindent\textbf{Proof. }%
$<$%
Case 1%
$>$%
. Assume that $E=+\infty$. From $\lim_{u\rightarrow+\infty}u/(u-\xi)=1$, for
any $0<\varepsilon<1/3$, there exists $N$ such that $1-\varepsilon
<u/(u-\xi)<1+\varepsilon$ for each $u>N.$ This implies that%
\[
\left\vert \frac{a(x)-u}{a(x)t-ut+u}\right\vert =\frac{1}{\left\vert
t-\frac{1}{1-\frac{a(x)}{u}}\right\vert }<\frac{1}{\varepsilon}\qquad
(\varepsilon\leq t\leq1-2\varepsilon,\text{ }u>N,\text{ }a(x)\neq u).
\]
Therefore, by Lebesgue (dominated-convergence) theorem, we obtain%
\[
\lim_{u\rightarrow+\infty}w_{u}(t)=\int_{I}\frac{1}{t-1}d(F(x))=\frac{-1}%
{1-t}<0.
\]
In particular, $\lim_{u\rightarrow+\infty}w_{u}(\varepsilon)=-1/(1-\varepsilon
)<0$. Therefore, there exists $M>0$ such that $w_{u}(\varepsilon
)<-1/(2(1-\varepsilon))$ for each $u>M$. On the basis of the fact that
$w_{u}(t)$ is strictly decreasing with respect to $t$, we have $0<\widetilde
{t}_{u}<\varepsilon$ for each $u$ $>M.$ This implies that $\lim_{u\rightarrow
+\infty}\widetilde{t}_{u}=0.$%

$<$%
Case 2%
$>$
Assume that $E<+\infty$. By Lemma 3.1, the analytic function $w_{E}(t)$ is
well defined with respect to $t\in(0,$ $E/(E-\xi))$. Similarly, from Lemmas
3.2 and 3.3, we know that $w_{E}(t)$ is strictly decreasing and $\lim
_{t\rightarrow0^{+}}w_{E}(t)=0.$ Therefore, we have $w_{E}(t)<0$.

If $0<\varepsilon<E/(2(E$ $-\xi))$ and $(E$ $+\max(0,\xi))/2<u<E$, then due to
$0<\varepsilon$ $<u/(u-\xi),$ $w_{u}(\varepsilon)$ is well defined. By Lemma
3.1 we have $\lim_{u\rightarrow E^{-}}w_{u}(\varepsilon)$ $=w_{E}%
(\varepsilon)<0$. Therefore, there exists $\delta>0$ such that $w_{u}%
(\varepsilon)<0$ for each $u\in(E-\delta,E)$. This implies that $0<\widetilde
{t}_{u}<\varepsilon$ and $\lim_{u\rightarrow E^{-}}\widetilde{t}_{u}=0.$%
\hfill$\square$

\bigskip

\noindent\textbf{Lemma 3.8.}\textit{ If }$\xi>0$\textit{, }$\widetilde
{t}_{1/H}=1.$

\bigskip

\noindent\textbf{Proof}. If $\xi>0$, then by Lemma 3.5 we have $\xi+1/H_{\xi
}<1/H<E$ and $1<u/(u-\xi)$. Therefore, $w_{u}(t)$ and $\widetilde{t}_{u}$\ are
analytic near $(u,t)=(1/H$, $1)$. The conclusion follows from the equality%
\[
w_{1/H}(1)=\int_{I}\frac{a(x)-\frac{1}{H}}{a(x)-\frac{1}{H}+\frac{1}{H}%
}d(F(x))=0.
\]

\hfill$\square$

\bigskip

\noindent\textbf{Lemma 3.9.} \textit{If }$\xi>0$\textit{, }$\lim
_{u\rightarrow(\xi+1/H_{\xi})^{+}}\widetilde{t}_{u}=1+\xi H_{\xi}.$

\bigskip

\noindent\textbf{Proof}. Due to Lemma 3.6, $\lim_{u\rightarrow(\xi+1/H_{\xi
})^{+}}\widetilde{t}_{u}$ exists, and we denote it by $\gamma$. It is clear
that $\gamma>1$. According to the inequality $\widetilde{t}_{u}<u/(u-\xi)$, we
have $\gamma\leq1+\xi H_{\xi}$.

Assume $H_{\xi}<+\infty$ and $\gamma<1+\xi H_{\xi}$, then, for any
$0<\varepsilon<\min((1+\xi H_{\xi}$ $-\gamma)/3,\gamma/2)$, there exists
$\delta>0$ such that%
\[
\left\vert \widetilde{t}_{u}-\gamma\right\vert <\varepsilon\text{ and
}\left\vert \frac{u}{u-\xi}-(1+\xi H_{\xi})\right\vert <\varepsilon
\]
for each $u\in(\xi+1/H_{\xi},\xi+1/H_{\xi}+\delta).$ This implies that
\[
\left\vert \widetilde{t}_{u}-\frac{1}{1-\frac{a(x)}{u}}\right\vert
>\varepsilon,
\]
where $a(x)\neq u$. By Lebesgue theorem, we obtain%
\[
0=\lim_{u\rightarrow(\xi+1/H_{\xi})^{+}}\int_{I}\frac{1}{\widetilde{t}%
_{u}-\frac{1}{1-\frac{a(x)}{u}}}d(F(x))=\int_{I}\frac{1}{\gamma-\frac
{1}{1-\frac{a(x)}{\xi+1/H_{\xi}}}}d(F(x)).
\]
This is a contradiction because the term on the right is positive, which is
deduced from the fact that the function%
\[
\int_{I}\frac{1}{t-\frac{1}{1-\frac{a(x)}{\xi+1/H_{\xi}}}}d(F(x))
\]
is strictly decreasing form $E/(\xi+1/H_{\xi})-1>0$ to $0$ with respect to
$t\in(0$, $1+\xi H_{\xi})$.

Assume $H_{\xi}=+\infty$ and $\gamma<+\infty$. Then, we have $0<1/(a(x)\gamma
-\xi\gamma+\xi)$ $\leq1/(a(x)\widetilde{t}_{u}-\xi\widetilde{t}_{u}+\xi
)\leq1/\xi$ for each $\xi<u<E$. Therefore, by Lebesgue theorem, we obtain%
\begin{align*}
0  &  =\lim_{u\rightarrow\xi^{+}}\int_{I}\frac{a(x)-u}{a(x)\widetilde{t}%
_{u}-u\widetilde{t}_{u}+u}d(F(x))\\
&  =\lim_{u\rightarrow\xi^{+}}\frac{1}{\widetilde{t}_{u}}(1-u\int_{I}\frac
{1}{a(x)\widetilde{t}_{u}-u\widetilde{t}_{u}+u}d(F(x)))\\
&  =\frac{1}{\gamma}(1-\xi\int_{I}\frac{1}{a(x)\gamma-\xi\gamma+\xi
}d(F(x)))>\frac{1}{\gamma}(1-\xi\int_{I}\frac{1}{\xi}d(F(x)))=0,
\end{align*}
which is a contradiction. This implies that if $H_{\xi}=+\infty$,
$\gamma=1+\xi H_{\xi}=+\infty.$

\hfill$\square$

\bigskip

\noindent\textbf{Lemma 3.10.} \textit{If }$\xi<0$\textit{ and }$\xi+1/H_{\xi
}>0,$\textit{ }$\lim_{u\rightarrow(\xi+1/H_{\xi})^{+}}\widetilde{t}_{u}=1+\xi
H_{\xi}.$

\bigskip

\noindent\textbf{Proof}. Due to $1/H_{\xi}>-\xi>0$, we have $H_{\xi}<+\infty$.
It should be noted that there exists $\delta>0$ such that $\widetilde{t}_{u}$
is strictly increasing or decreasing in the interval $u\in(\xi+1/H_{\xi}$,
$\xi+1/H_{\xi}+\delta)$, which is demonstrated in the proof of Lemma 3.16.
Therefore, $\lim_{u\rightarrow(\xi+1/H_{\xi})^{+}}\widetilde{t}_{u}$ exists
and denoted by $\gamma$. By the inequality $\widetilde{t}_{u}<u/(u-\xi)$, we
have $\gamma\leq1+\xi H_{\xi}$. Assume $\gamma<$ $1+\xi H_{\xi}$, then, as in
the proof of Lemma 3.9, we have a contradiction.\hfill$\square$

\bigskip

\noindent\textbf{Lemma 3.11.} \textit{If }$\xi=0$\textit{ and }$1/H>0,$%
\textit{ }$\lim_{u\rightarrow(1/H)^{+}}\widetilde{t}_{u}=1.$

\bigskip

\noindent\textbf{Proof}. It should be noted that $H<+\infty$. Due to Lemma
3.6, $\lim_{u\rightarrow(1/H)^{+}}\widetilde{t}_{u}$ exists, and is denoted by
$\gamma$. According to the relation $\widetilde{t}_{u}<u/(u-\xi)=1$, we have
$\gamma\leq1$. Assume $\gamma<$ $1$, then as in the proof of Lemma 3.9, the
function mentioned there is strictly decreasing from $HE-1>0$ to $0$ in the
interval $t\in(0$, $1),$ which leads to a contradiction.\hfill$\square$

\bigskip

\noindent\textbf{Lemma 3.12.} \textit{If }$\xi<0$\textit{ and }$\xi+1/H_{\xi
}\leq0,$\textit{ }$\lim_{u\rightarrow0^{+}}\widetilde{t}_{u}=0.$

\bigskip

\noindent\textbf{Proof}. On the basis of the definition, $0<\widetilde{t}%
_{u}<u/(u-\xi)$ and $\max(0$, $\xi)=0<u$ $<E$. Therefore, $0\leq
\lim_{u\rightarrow0^{+}}\widetilde{t}_{u}\leq\lim_{u\rightarrow0^{+}}%
u/(u-\xi)=0.$\hfill$\square$

\bigskip

\noindent\textbf{Lemma 3.13.} \textit{If }$\xi\leq0$\textit{, }$\widetilde
{t}_{u}<\int_{a(x)\neq0}d(F(x)).$

\bigskip

\noindent\textbf{Proof. }As $\xi\leq0$, we have $0<\widetilde{t}_{u}%
<u/(u-\xi)\leq1$ for each $u\in(\xi+1/H_{\xi}$, $E)$. In this case, the
equation $w_{u}(\widetilde{t}_{u})=0$ is equivalent to
\begin{align*}
&  \frac{-1}{-\widetilde{t}_{u}+1}\int_{a(x)=0}d(F(x))+\frac{1}{\widetilde
{t}_{u}}\int_{a(x)\neq0}d(F(x))\\
&  =\frac{u}{\widetilde{t}_{u}}\int_{a(x)\neq0}\frac{1}{a(x)\widetilde{t}%
_{u}-u\widetilde{t}_{u}+u}d(F(x))>0,
\end{align*}
which leads to $\widetilde{t}_{u}<\int_{a(x)\neq0}d(F(x)).$\hfill$\square$

\bigskip

\noindent\textbf{Lemma 3.14.}\textit{ If }$\xi=0$\textit{ and }$H=+\infty
$\textit{, }$\lim_{u\rightarrow0^{+}}\widetilde{t}_{u}=\int_{a(x)>0}%
d(F(x))$\textit{.}

\bigskip

\noindent\textbf{Proof.} Due to Lemma 3.6, $\lim_{u\rightarrow0^{+}}%
\widetilde{t}_{u}$ exists, and we denote it by $\gamma$. We can choose
$\delta>0$ such that $\gamma/2<\widetilde{t}_{u}<\gamma$ for each $u\in(0,$
$\delta)$. It should be noted that the equation $w_{u}(\widetilde{t}_{u})=0$
is equivalent to%
\[
\frac{-1}{-\widetilde{t}_{u}+1}\int_{a(x)=0}d(F(x))+\int_{a(x)>0}\frac
{1}{\widetilde{t}_{u}-\frac{1}{1-\frac{a(x)}{u}}}d(F(x))=0.
\]

Assume that $\gamma<1$, we have
\[
\left\vert \frac{1}{\widetilde{t}_{u}-\frac{1}{1-\frac{a(x)}{u}}}\right\vert
<\left\{
\begin{array}
[c]{cc}%
\frac{1}{1-\gamma}, & \text{if }0<a(x)<u,\\
\frac{2}{\gamma}, & \text{if }a(x)>u,
\end{array}
\right.
\]
for each $u\in(0,$ $\delta)$. In this case, by Lebesgue theorem, we obtain%
\[
\frac{-1}{-\gamma+1}\int_{a(x)=0}d(F(x))+\frac{1}{\gamma}\int_{a(x)>0}%
d(F(x))=0,
\]
which implies that $\gamma=\int_{a(x)>0}d(F(x)).$

In the case in which$\ \gamma=1$, due to Lemma 3.13, we have $\gamma\leq
\int_{a(x)>0}d(F(x))$ $\leq1$. Thus, $\gamma=\int_{a(x)>0}d(F(x)).$%
\hfill$\square$

\bigskip

\noindent\textbf{Lemma 3.15. }\textit{The function }$\widetilde{t}_{u}%
/u$\textit{ is strictly decreasing with respect to }$u\in(\max(0$\textit{,
}$\xi+1/H_{\xi}),$\textit{ }$E).$

\bigskip

\noindent\textbf{Proof.} Using the equality%
\[
\int_{I}\frac{a(x)-u}{a(x)\widetilde{t}_{u}-u\widetilde{t}_{u}+u}%
d(F(x))=\int_{I}\frac{(a(x)-u)^{2}\widetilde{t}_{u}+a(x)u-u^{2}}{\left(
a(x)\widetilde{t}_{u}-u\widetilde{t}_{u}+u\right)  ^{2}}d(F(x))=0
\]
and by the proof of Lemma 3.6, we obtain%
\begin{align*}
\frac{d}{du}\left(  \frac{\widetilde{t}_{u}}{u}\right)   &  =-\frac{u\int
_{I}\frac{a(x)}{(a(x)\widetilde{t}_{u}-u\widetilde{t}_{u}+u)^{2}%
}d(F(x))+\widetilde{t}_{u}\int_{I}\frac{(a(x)-u)^{2}}{(a(x)\widetilde{t}%
_{u}-u\widetilde{t}_{u}+u)^{2}}d(F(x))}{u^{2}\int_{I}\frac{(a(x)-u)^{2}%
}{(a(x)\widetilde{t}_{u}-u\widetilde{t}_{u}+u)^{2}}d(F(x))}\\
&  =-\frac{\int_{I}\frac{1}{(a(x)\widetilde{t}_{u}-u\widetilde{t}_{u}+u)^{2}%
}d(F(x))}{\int_{I}\frac{(a(x)-u)^{2}}{(a(x)\widetilde{t}_{u}-u\widetilde
{t}_{u}+u)^{2}}d(F(x))}<0.
\end{align*}

\hfill$\square$

\bigskip

We define the continuous function $\overline{t}_{u}$ in the interval $[0$,
$+\infty)$ as follows:%
\begin{equation}
\text{If }\xi>0,\text{ then }\overline{t}_{u}:=\left\{
\begin{array}
[c]{cc}%
1, & \text{if }0\leq u\leq1/H,\\
\widetilde{t}_{u}, & \text{if }1/H<u<E,\\
0, & \text{if }u\geq E.
\end{array}
\right.  \tag{7}%
\end{equation}%
\begin{equation}
\text{If }\xi=0,\text{ then }\overline{t}_{u}:=\left\{
\begin{array}
[c]{cc}%
\int_{a(x)>0}d(F(x)), & \text{if }u=0.\\
1, & \text{if }0<u\leq1/H,\\
\widetilde{t}_{u}, & \text{if }1/H<u<E,\\
0, & \text{if }u\geq E.
\end{array}
\right.  \tag{8}%
\end{equation}%
\begin{equation}
\text{If }\xi<0,\text{ then }\overline{t}_{u}:=\left\{
\begin{array}
[c]{cc}%
\frac{u}{u-\xi}, & \text{if }0\leq u\leq\max(0,\text{ }\xi+1/H_{\xi}),\\
\widetilde{t}_{u}, & \text{if }\max(0,\text{ }\xi+1/H_{\xi})<u<E,\\
0, & \text{if }u\geq E.
\end{array}
\right.  \tag{9}%
\end{equation}
\bigskip

\noindent\textbf{Lemma 3.16}. \textit{If} $\xi<0$, \textit{then the value}
$0<u_{\max}<E$ \textit{exists, which satisfies the following properties:}

(a) $\overline{t}_{u}$ \textit{is strictly increasing in the interval}
$0<u<u_{\max}.$

(b) $\overline{t}_{u}$ \textit{is strictly decreasing in the interval}
$u_{\max}<u<E.$

\bigskip

\noindent\textbf{Proof}. $\overline{t}_{u}/u=1/(u-\xi)$ is strictly decreasing
in the interval $0<u<\xi+1/H_{\xi}$, if $\xi+1/H_{\xi}>0$. Using Lemma 3.15,
we have that $\overline{t}_{u}/u$ is strictly decreasing in the interval
$0<u<E$. We denote the value of $\lim_{u\rightarrow0^{+}}\overline{t}_{u}/u$
by $\eta>0.$ As the function $y:=\overline{t}_{u}/u$ is strictly decreasing
with respect to $u$, $\overline{t}_{u}$ can be considered to be a function
with a variable $y\in(0$, $\eta)$.

If $u\in(\max(0$, $\xi+1/H_{\xi})$, $E)$, $\overline{t}_{u}=\widetilde{t}_{u}%
$. From $w_{u}(\widetilde{t}_{u})=0$ we have%
\[
\int_{I}\frac{1}{a(x)y-\overline{t}_{u}+1}d(F(x))=1.
\]
Thus,
\[
\int_{I}\frac{a(x)-\frac{d\overline{t}_{u}}{dy}}{\left(  a(x)y-\overline
{t}_{u}+1\right)  ^{2}}d(F(x))=0.
\]
This implies that
\[
\frac{d\overline{t}_{u}}{dy}=\frac{\int_{I}\frac{a(x)}{(a(x)y-\overline{t}%
_{u}+1)^{2}}d(F(x))}{\int_{I}\frac{1}{(a(x)y-\overline{t}_{u}+1)^{2}}%
d(F(x))}.
\]
Denoting $\int_{I}\frac{1}{(a(x)y-\overline{t}_{u}+1)^{2}}d(F(x))$ by $s$, we
obtain%
\begin{align*}
\frac{d^{2}\overline{t}_{u}}{dy^{2}}  &  =\frac{1}{s^{2}}\left(
\begin{array}
[c]{c}%
-2s\int_{I}\frac{a(x)(a(x)-\frac{d\overline{t}_{u}}{dy})}{(a(x)y-\overline
{t}_{u}+1)^{3}}d(F(x))\\
+2\int_{I}\frac{a(x)}{(a(x)y-\overline{t}_{u}+1)^{2}}d(F(x))\times\int
_{I}\frac{a(x)-\frac{d\overline{t}_{u}}{dy}}{(a(x)y-\overline{t}_{u}+1)^{3}%
}d(F(x))
\end{array}
\right) \\
&  =\frac{-2}{s^{3}}\left(
\begin{array}
[c]{c}%
s^{2}\int_{I}\frac{a(x)^{2}}{(a(x)y-\overline{t}_{u}+1)^{3}}d(F(x))\\
-2s\int_{I}\frac{a(x)}{(a(x)y-\overline{t}_{u}+1)^{2}}d(F(x))\int_{I}%
\frac{a(x)}{(a(x)y-\overline{t}_{u}+1)^{3}}d(F(x))\\
+\left(  \int_{I}\frac{a(x)}{(a(x)y-\overline{t}_{u}+1)^{2}}d(F(x))\right)
^{2}\times\int_{I}\frac{1}{(a(x)y-\overline{t}_{u}+1)^{3}}d(F(x))
\end{array}
\right)  .
\end{align*}
As a quadratic function with respect to $s$, $-s_{3}/2\times d^{2}\overline
{t}_{u}/dy^{2}$ has the determinant given by%
\begin{align*}
&  \left(  \int_{I}\frac{a(x)}{(a(x)y-\overline{t}_{u}+1)^{2}}d(F(x))\right)
^{2}\times\\
&  \left(
\begin{array}
[c]{c}%
\left(  \int_{I}\frac{a(x)}{(a(x)y-\overline{t}_{u}+1)^{3}}d(F(x))\right)
^{2}\\
-\int_{I}\frac{a(x)^{2}}{(a(x)y-\overline{t}_{u}+1)^{3}}d(F(x))\times\int
_{I}\frac{1}{(a(x)y-\overline{t}_{u}+1)^{3}}d(F(x))
\end{array}
\right)
\end{align*}
Due to H\"{o}lder inequality with respect to the two functions $\frac
{a(x)}{(a(x)y-\overline{t}_{u}+1)^{3/2}}$ and $\frac{1}{(a(x)y-\overline
{t}_{u}+1)^{3/2}}$, this determinant is negative. Therefore, we have
$\frac{d^{2}\overline{t}_{u}}{dy^{2}}<0$, which implies that $\frac
{d\overline{t}_{u}}{dy}$ is strictly decreasing.

First, we assume that $\xi+1/H_{\xi}\leq0$. Assign $\alpha:=\lim
_{y\rightarrow0^{+}}d\overline{t}_{u}/dy$ and $\beta$ $:=\lim_{y\rightarrow
\eta^{-}}d\overline{t}_{u}/dy.$

If $\alpha\leq0$, then $d\overline{t}_{u}/dy<0$ for each $0<y<\eta$. This
contradicts the fact that $\overline{t}_{u}=0$ if $y=0$ (Lemma 3.7), and
$\overline{t}_{u}>0$ if $0<y<\eta$. Therefore, $\alpha>0.$

If $\beta\geq0,$ $d\overline{t}_{u}/dy>0$ for each $0<y<\eta$. This
contradicts the fact that $\overline{t}_{u}=0$ if $y=\eta$ (Lemma 3.11), and
$\overline{t}_{u}>0$ if $0<y<\eta$. Therefore, $\beta<0.$

The value $0<y_{\max}<\eta$ such that $d\overline{t}_{u}/dy|_{y=y_{\max}}=0$
can then be determined. The value $0<u_{\max}<E$ required is determined using
$y_{\max}.$

Second, we assume that $\xi+1/H_{\xi}>0$. If $u\in(0,$ $\xi+1/H_{\xi})$, then
$\overline{t}_{u}$ $=u/(u-\xi)$ is strictly increasing. Moreover, we have
$u/(u-\xi)|_{u=\xi+1/H_{\xi}}=1+\xi H_{\xi}=\lim_{u\rightarrow(\xi+1/H_{\xi
})^{+}}\widetilde{t}_{u}$ (Lemma 3.10) and $y|_{u=\xi+1/H_{\xi}}=H_{\xi}$.
Redefine $\beta$ as $\lim_{y\rightarrow H_{\xi}^{-}}d\overline{t}_{u}/dy$. If
$\beta<0$, then as above, we obtain the required value $H_{\xi}$%
$<$%
$u_{\max}$%
$<$%
$E$. If $\beta\geq0$, $d\overline{t}_{u}/dy>0$ for each $0<y<H_{\xi}$. This
implies that $d\overline{t}_{u}/du$ $=d\overline{t}_{u}/dy\times dy/du<0$
since $dy/du<0$ for each $u\in(\xi+1/H_{\xi},$ $E)$ (Lemma 3.15). Thus,
$\overline{t}_{u}$ is strictly decreasing. Therefore, the required value is
$u_{\max}=\xi+1/H_{\xi}$.\hfill$\square$

\bigskip

\noindent{\large \textbf{4. Pre-growth rate}}

\bigskip

In this Section we assume that $u\in(\max(0$, $\xi),$ $E)$ and $\rho,$
$t\in(0,$ $u/(u-\xi))$ unless otherwise mentioned. Define $G_{u,\rho}(t)$ by
the equality%
\begin{equation}
\exp\left(  \int_{\rho}^{t}w_{u}(t)dt\right)  =\exp\left(  \int_{I}\log
\frac{a(x)t-ut+u}{a(x)\rho-u\rho+u}d(F(x))\right)  , \tag{10}%
\end{equation}
which can be verified using the following inequalities.%
\[
\min\left(  \frac{t}{\rho},\text{ }\frac{u/(u-\xi)-t}{u/(u-\xi)-\rho}\right)
<\left\vert \frac{a(x)t-ut+u}{a(x)\rho-u\rho+u}\right\vert <\max\left(
\frac{t}{\rho},\text{ }\frac{u/(u-\xi)-t}{u/(u-\xi)-\rho}\right)  ,
\]%
\[
\left\vert \frac{a(x)-u}{a(x)s-us+u}\right\vert <\frac{1}{\min(\rho,\text{
}t,\text{ }u/(u-\xi)-\rho,\text{ }u/(u-\xi)-t)}%
\]
for each $x\in I$ and $s\in(\min(\rho$, $t),$ $\max(\rho$, $t))$. As
$w_{u}(t)$ is strictly decreasing with respect to $t$ from the positive value
$E/u-1,$ to the value $(1-\xi/u)H_{\xi}(\xi+1/H_{\xi}-u)$ (Lemmas 3.2, 3.3,
and 3.4), $\int_{\rho}^{t}w_{u}(t)dt$ is strictly decreasing with respect to
$\rho$ near $0^{+}$. Therefore, the limit%
\begin{align}
\lim_{\rho\rightarrow0^{+}}\exp\left(  \int_{\rho}^{t}w_{u}(t)dt\right)   &
=\exp\left(  \lim_{\rho\rightarrow0^{+}}\int_{I}\log\frac{a(x)t-ut+u}%
{a(x)\rho-u\rho+u}d(F(x))\right) \tag{11}\\
&  =\exp\left(  \int_{I}\log\left(  a(x)t/u-t+1\right)  d(F(x))\right)
\nonumber
\end{align}
finitely exists or $+\infty$, which we denote by $\widetilde{G}_{u}(t)$ and
refer to as \textit{pre-growth rate}. The equality mentioned above is obtained
using Lebesgue theorem because the integrand is monotone with respect to
$\rho$ in $\{$ $x$ $|$ $a(x)>u$ $\}$ or $\{$ $x$ $|$ $a(x)$%
$<$%
$u$ $\}$.

\bigskip

\noindent\textbf{Lemma 4.1.} $\widetilde{G}_{u}(t)<E/u.$

\bigskip

\noindent\textbf{Proof}. By Jensen's inequality we have%
\begin{align*}
\int_{I}\log\left(  a(x)t/u-t+1\right)  d(F(x))  &  \leq\log\int_{I}\left(
a(x)t/u-t+1\right)  d(F(x))\\
&  =\log(Et/u-t+1)<\log(E/u).
\end{align*}

\hfill$\square$

\bigskip

\noindent\textbf{Lemma 4.2.} $\int_{a(x)<u}\left\vert \log\left(
a(x)t/u-t+1\right)  \right\vert d(F(x))<+\infty.$

\bigskip

\noindent\textbf{Proof}. In general, $a(x)t/u-t+1\geq(u-\xi)(u/(u-\xi
)-t)/u>0$. If $a(x)<u$, $a(x)t/u-t+1$%
$<$%
$1$. Therefore, we obtain%
\[
\int_{a(x)<u}\left\vert \log\left(  a(x)t/u-t+1\right)  \right\vert
d(F(x))\leq\left\vert \log\left(  \frac{u-\xi}{u}\left(  \frac{u}{u-\xi
}-t\right)  \right)  \right\vert <+\infty.
\]

\hfill$\square$

\bigskip

\noindent\textbf{Lemma 4.3.} \textit{The following three statements are
equivalent.}

(1)\textit{ }$\int_{a(x)>1}\log a(x)d(F(x))<+\infty.$

(2)\textit{ }$\widetilde{G}_{u}(t)<+\infty$\textit{ for each }$u$\textit{ and
}$t.$

(3)\textit{ }$\widetilde{G}_{u_{1}}(t_{1})<+\infty$\textit{ for some }$u_{1}%
$\textit{ and }$t_{1}.$

\bigskip

\noindent\textbf{Proof}. (1) $\Longrightarrow$ (2). The function $\log\left(
a(x)t/u-t+1\right)  $ satisfies the following inequalities.%
\[
\int_{u<a(x)\leq1\text{ or }1\leq a(x)<u}\left\vert \log a(x)\right\vert
d(F(x))\leq\left\vert \log u\right\vert <+\infty.
\]%
\begin{align*}
&  \left\vert \int_{a(x)>u}\log\left(  a(x)t/u-t+1\right)  d(F(x))-\int
_{a(x)>u}\log a(x)d(F(x))\right\vert \\
&  =\left\vert \int_{a(x)>u}\left(  \log\frac{t}{u}+\log\left(  1+\frac
{u(1-t)}{a(x)t}\right)  \right)  d(F(x))\right\vert \\
&  \leq\left\vert \log\frac{t}{u}\right\vert +\left\vert \log t\right\vert
<+\infty.
\end{align*}
Based on Lemma 4.2, we obtain the integrability of $\log\left(
a(x)t/u-t+1\right)  .$

(3) $\Longrightarrow$ (1). It should be noted that $u_{1}\in(\max(0$, $\xi),$
$E)$ and $t_{1}\in(0,$ $u_{1}/(u_{1}-\xi)).$ The result can be obtained in a
similar manner as above.

(2) $\Longrightarrow$ (3). It is clear.\hfill$\square$

\bigskip

If one of the above three statements is satisfied, we can write $\widetilde
{G}<+\infty$.

\bigskip

\noindent\textbf{Lemma 4.4.} \textit{If} $\widetilde{G}<+\infty,$
$\lim_{t\rightarrow0^{+}}\widetilde{G}_{u}(t)=1$.

\bigskip

\noindent\textbf{Proof. }Since $\lim_{t\rightarrow0^{+}}w_{u}(t)=E/u-1>0$
(Lemma 3.3), $\int_{0}^{t}w_{u}(t)dt$ is strictly increasing and bounded with
respect to $t$ near $0^{+}$. Therefore, we obtain that $\lim_{t\rightarrow
0^{+}}$ $\int_{0}^{t}w_{u}(t)dt=0.$\hfill$\square$

\bigskip

\noindent\textbf{Lemma 4.5.} \textit{If }$u\in(\max(0,$\textit{ }$\xi
+1/H_{\xi})$, $E)$\textit{, }$\max_{0<t<u/(u-\xi)}G_{u,\rho}(t)=G_{u,\rho
}(\widetilde{t}_{u})$.

\bigskip

\noindent\textbf{Proof. }It is clear form the facts that\textbf{ }%
$0<G_{u,\rho}(t)<+\infty$ and\textbf{ }%
\begin{align*}
\partial G_{u,\rho}(t)/\partial t  &  =\frac{\partial}{\partial t}\exp\left(
\int_{\rho}^{t}w_{u}(t)dt\right) \\
&  =G_{u,\rho}(t)\left(  \frac{\partial}{\partial t}\int_{\rho}^{t}%
w_{u}(t)dt\right)  =G_{u,\rho}(t)w_{u}(t).
\end{align*}

\hfill$\square$

\bigskip

\noindent\textbf{Lemma 4.6.} \textit{If }$\widetilde{G}<+\infty$\textit{ and
}$u\in(\max(0,$\textit{ }$\xi+1/H_{\xi})$, $E)$\textit{, }$\max_{0<t<u/(u-\xi
)}$\textit{ }$\widetilde{G}_{u}(t)=\widetilde{G}_{u}(\widetilde{t}_{u}).$

\bigskip

\noindent\textbf{Proof. }In a similar manner as that of the proof of Lemma
4.5, we have
\begin{align*}
\frac{\partial\widetilde{G}_{u}(t)}{\partial t}  &  =\frac{\partial}{\partial
t}\lim_{\rho\rightarrow0^{+}}\exp\left(  \int_{\rho}^{t}w_{u}(t)dt\right)
=\lim_{\rho\rightarrow0^{+}}\frac{\partial}{\partial t}\exp\left(  \int_{\rho
}^{t}w_{u}(t)dt\right) \\
&  =\lim_{\rho\rightarrow0^{+}}G_{u,\rho}(t)w_{u}(t)=\widetilde{G}_{u}%
(t)w_{u}(t),
\end{align*}
which implies the conclusion.\hfill$\square$

\bigskip

\noindent\textbf{Lemma 4.7.} \textit{Two functions }$G_{u,\rho}(t)$\textit{
and }$\widetilde{G}_{u}(t)$\textit{ }$(<+\infty)$\textit{ are concave with
respect to }$t$\textit{.}

\bigskip

\noindent\textbf{Proof. }Using Lemmas 3.2, 4.5, and H\"{o}lder inequality, we
have%
\begin{align*}
\partial^{2}G_{u,\rho}(t)/\partial t^{2}  &  =G_{u,\rho}(t)(w_{u}%
^{2}(t)+\partial w_{u}(t)/\partial t)\\
&  =G_{u,\rho}(t)\left(
\begin{array}
[c]{c}%
\left(  \int_{I}(a(x)-u)/(a(x)t-ut+u)d(F(x))\right)  ^{2}\\
-\int_{I}(a(x)-u)^{2}/(a(x)t-ut+u)^{2}d(F(x))
\end{array}
\right)  <0.
\end{align*}

\hfill$\square$

\bigskip

Using the above result, we also have $\partial^{2}\widetilde{G}_{u}%
(t)/\partial t^{2}<0.$

\bigskip

\noindent\textbf{Lemma 4.8.} \textit{If }$\xi\geq0$ \textit{and} $t>\rho
$\textit{, }$G_{u,\rho}(t)$\textit{ is strictly decreasing with respect to
}$u$\textit{.}

\bigskip

\noindent\textbf{Proof. }From $a(x)\geq0$, we obtain%
\begin{align*}
\frac{\partial G_{u,\rho}(t)}{\partial u}  &  =\frac{\partial}{\partial u}%
\exp\left(  \int_{\rho}^{t}w_{u}(t)dt\right)  =G_{u,\rho}(t)\times
\frac{\partial}{\partial u}\int_{\rho}^{t}w_{u}(t)dt\\
&  =-G_{u,\rho}(t)\times\int_{\rho}^{t}\left(  \int_{I}\frac{a(x)}%
{(a(x)t-ut+u)^{2}}d(F(x))\right)  dt<0.
\end{align*}

\hfill$\square$

\bigskip

\noindent\textbf{Lemma 4.9.}\textit{ If }$\xi\geq0$\textit{ and }%
$\widetilde{G}<+\infty$\textit{, }$\widetilde{G}_{u}(t)$\textit{ is strictly
decreasing with respect to }$u$\textit{.}

\bigskip

\noindent\textbf{Proof. }Using Lemma 4.8, we obtain the conclusion.\hfill
$\square$

\bigskip

\noindent\textbf{Lemma 4.10. }\textit{If }$\widetilde{G}<+\infty$\textit{,
}$\lim_{t\rightarrow(u/(u-\xi))^{-}}\widetilde{G}_{u}(t)=\exp(\int_{I}%
\log\left(  a(x)-\xi\right)  d(F(x)))$\textit{ }$/(u-\xi)$\textit{ for each
}$u.$

\bigskip

\noindent\textbf{Proof. }If $a(x)>u$, $a(x)t/u-t+1$ is strictly increasing
with respect to $t$. Therefore, using Lebesgue theorem, we have%
\begin{align*}
0  &  \leq\lim_{t\rightarrow(u/(u-\xi))^{-}}\int_{a(x)>u}\log
(a(x)t/u-t+1)d(F(x))\\
&  =\int_{a(x)>u}\log\frac{a(x)-\xi}{u-\xi}d(F(x))<+\infty.
\end{align*}
On the other hand, if $a(x)<u$, then $a(x)t/u-t+1$ is strictly decreasing with
respect to $t$. Hence, using Lebesgue theorem, we have%
\begin{align*}
&  \lim_{t\rightarrow(u/(u-\xi))^{-}}\int_{a(x)<u}\log(a(x)t/u-t+1)d(F(x))\\
&  =\int_{a(x)<u}\log\frac{a(x)-\xi}{u-\xi}d(F(x))<0,
\end{align*}
which implies that%
\begin{align*}
&  \lim_{t\rightarrow(u/(u-\xi))^{-}}\widetilde{G}_{u}(t)\\
&  =\exp\left(  \int_{a(x)>u}\log\frac{a(x)-\xi}{u-\xi}d(F(x))+\int
_{a(x)<u}\log\frac{a(x)-\xi}{u-\xi}d(F(x))\right) \\
&  =\exp\left(  \int_{I}\log\left(  a(x)-\xi\right)  d(F(x))\right)  /(u-\xi).
\end{align*}

\hfill$\square$

\bigskip

As an expansion of the definition of $\widetilde{G}_{u}(t)$, we define
$\widetilde{G}_{u}((u/(u-\xi))^{-})$ by $\exp(\int_{I}\log\left(
a(x)-\xi\right)  d(F(x)))$ $/(u-\xi)$ for each $u\in(\max(0,$ $\xi)$, $E)$.

\bigskip

\noindent\textbf{Lemma 4.11.}\textit{ }$\widetilde{G}_{u}(\widetilde{t}%
_{u})>1$\textit{ if }$u\in(\max(0$\textit{, }$\xi+1/H_{\xi})$\textit{, }$E).$

\bigskip

\noindent\textbf{Proof.} If $0<t<\widetilde{t}_{u}$, $w_{u}(t)>0$. Hence, we
have
\[
\widetilde{G}_{u}(\widetilde{t}_{u})=\exp(\int_{0}^{\widetilde{t}_{u}}%
w_{u}(t)dt)>e^{0}=1.
\]

\hfill$\square$

\bigskip

\noindent\textbf{Lemma 4.12.} $\widetilde{G}_{u}(\widetilde{t}_{u})$\textit{
}$(<+\infty)$\textit{\ is strictly decreasing with respect to }$u\in(\max
(0$\textit{, }$\xi+1/H_{\xi})$\textit{, }$E).$

\bigskip

\noindent\textbf{Proof.} If $\left\vert a(x)\right\vert \leq2\left\vert
\xi\right\vert $, we have
\[
\left\vert \frac{a(x)}{a(x)t-ut+u}\right\vert =\left\vert \frac{a(x)}%
{(a(x)-\xi)t+\xi t-ut+u}\right\vert \leq\frac{2\left\vert \xi\right\vert
}{(u-\xi)(u/(u-\xi)-t)}.
\]
On the other hand, if $a(x)>2\left\vert \xi\right\vert $, we have
\[
\left\vert \frac{a(x)}{a(x)t-ut+u}\right\vert =\left\vert \frac{a(x)}%
{(a(x)-\xi)t+\xi t-ut+u}\right\vert \leq\left\vert \frac{a(x)}{(a(x)-\xi
)t}\right\vert <\frac{2}{t}.
\]
Thus, by the definition of $\widetilde{G}_{u}(t)$, we have%
\begin{align*}
\frac{\partial\widetilde{G}_{u}(t)}{\partial u}  &  =\widetilde{G}_{u}%
(t)\int_{I}\frac{\partial}{\partial u}\log(a(x)t/u-t+1)d(F(x))\\
&  =-\frac{t}{u}\widetilde{G}_{u}(t)\times\left(  \int_{I}\frac{a(x)}%
{a(x)t-ut+u}d(F(x))\right)  .
\end{align*}
The definition $w_{u}(\widetilde{t}_{u})=0$ leads to $\int_{I}%
a(x)/(a(x)\widetilde{t}_{u}-u\widetilde{t}_{u}+u)d(F(x))=1$. Therefore,%
\[
\frac{\partial\widetilde{G}_{u}(\widetilde{t}_{u})}{\partial u}=-\frac
{\widetilde{t}_{u}}{u}\widetilde{G}_{u}(\widetilde{t}_{u})<0.
\]

\hfill$\square$

\bigskip

\noindent\textbf{Lemma 4.13.} \textit{If }$\widetilde{G}<+\infty$\textit{,
}$\lim_{u\rightarrow E^{-}}\widetilde{G}_{u}(\widetilde{t}_{u})=1.$

\bigskip

\noindent\textbf{Proof.} From Lemmas 4.11 and 4.12, $\lim_{u\rightarrow E^{-}%
}\widetilde{G}_{u}(\widetilde{t}_{u})\geq1$ exists. Assume that $\xi\geq0$,
then from Lemmas 3.6 and 3.7, $\widetilde{t}_{u}/(1-\widetilde{t}_{u})$ is
strictly decreasing near $u=E^{-}$. Applying Lebesgue theorem to the equality%
\[
\widetilde{G}_{u}(\widetilde{t}_{u})-1+\widetilde{t}_{u}=(1-\widetilde{t}%
_{u})\left(  \exp\left(  \int_{I}\log\left(  \frac{a(x)\widetilde{t}_{u}%
}{u(1-\widetilde{t}_{u})}+1\right)  d(F(x))\right)  -1\right)
\]
we have
\[
\lim_{u\rightarrow E^{-}}\left(  \widetilde{G}_{u}(\widetilde{t}%
_{u})-1+\widetilde{t}_{u}\right)  =\exp\left(  \int_{I}\log\left(  0+1\right)
d(F(x))\right)  -1=0.
\]
This implies that $\lim_{u\rightarrow E^{-}}\widetilde{G}_{u}(\widetilde
{t}_{u})=1$.

In the case in which $\xi<0$, Lemma 3.16 can be used as a substitution of
Lemma 3.6 near $u=E^{-}$, where $\widetilde{t}_{u}$\ is strictly decreasing to
$0$. In order to apply Lebesgue theorem, it is sufficient to divide the above
integration into two parts $\{$ $x$ $|$ $a(x)\geq0$ $\}$ and $\{$ $x$ $|$
$a(x)<0$ $\}$.\hfill$\square$

\bigskip

\noindent\textbf{Lemma 4.14.} \textit{If }$\widetilde{G}<+\infty$\textit{ and
}$u\in(\max(0$\textit{, }$\xi+1/H_{\xi})$\textit{, }$\widetilde{G}%
_{u}(\widetilde{t}_{u})$ $=\exp\left(  \int_{u}^{E}\widetilde{t}%
_{u}/udu\right)  .$

\bigskip

\noindent\textbf{Proof}. Using $\partial\widetilde{G}_{u}(\widetilde{t}%
_{u})/\partial u=-\widetilde{t}_{u}/u\times\widetilde{G}_{u}(\widetilde{t}%
_{u})$\ (Lemmas 4.12 and 4.13), we can solve the differential equation.\hfill
$\square$

\bigskip

\noindent\textbf{Lemma 4.15} \textit{If }$\xi=0$\textit{ and }$H=+\infty
$\textit{, }$\lim_{u\rightarrow0^{+}}\widetilde{G}_{u}(\widetilde{t}%
_{u})=+\infty.$

\bigskip

\noindent\textbf{Proof}. Lemma 4.12 ensures the existence of $\lim
_{u\rightarrow0^{+}}\widetilde{G}_{u}(\widetilde{t}_{u})$, which is finite or
$+\infty$. If $a(x)>0$, $a(x)/u$ is strictly decreasing with respect to $u$.
Using Lebesgue theorem, we have
\begin{align*}
\lim_{u\rightarrow0^{+}}\widetilde{G}_{u}(\widetilde{t}_{u})  &
\geq\underline{\lim}_{u\rightarrow0^{+}}\widetilde{G}_{u}(\frac{1}{2})\\
&  \geq\underline{\lim}_{u\rightarrow0^{+}}\frac{1}{2}\exp\left(
\int_{a(x)>0}\log\frac{a(x)+u}{2u}d(F(x))\right)  =+\infty.
\end{align*}

\hfill$\square$

\bigskip

\noindent\textbf{Lemma 4.16}.\textit{ If }$\xi=0$\textit{, }$H<+\infty
,$\textit{ and }$\widetilde{G}<+\infty$\textit{, }$\lim_{u\rightarrow
(1/H)^{+}}\widetilde{G}_{u}(\widetilde{t}_{u})$ $=H\exp\left(  \int_{I}\log
a(x)d(F(x))\right)  .$

\bigskip

\noindent\textbf{Proof.} By definition, $\widetilde{G}<+\infty$ implies that
$\int_{a(x)>1}$ $\log a(x)d(F(x))<+\infty$. From Jensen's inequality theorem,
we have%
\[
+\infty>\log H=\log\int_{I}\frac{1}{a(x)}d(F(x))\geq\int_{I}\log\frac{1}%
{a(x)}d(F(x)),
\]
which implies that $\int_{I}\log a(x)d(F(x))>-\infty$. Therefore, $\log a(x)$
is integrable.

It should be noted that $\lim_{u\rightarrow(1/H)^{+}}\widetilde{t}_{u}=1$
(Lemma 3.11). Using the equalities $\lim_{u\rightarrow(1/H)^{+}}$
$\widetilde{t}_{u}/u=H$ and $\lim_{u\rightarrow(1/H)^{+}}(1-\widetilde{t}%
_{u})u/\widetilde{t}_{u}=0,$ we can choose $0<\delta<\min(1/H,$ $E-1/H)$, such
that $H/2<\widetilde{t}_{u}/u<3H/2$ and $(1-\widetilde{t}_{u})u/\widetilde
{t}_{u}<1/2$ for each $u\in(1/H,$ $1/H+\delta)$. Therefore, we have the
following properties.

(1) If $a(x)\geq1/2$, then%
\begin{align*}
\left\vert \log\left(  a(x)\widetilde{t}_{u}/u-\widetilde{t}_{u}+1\right)
\right\vert  &  =\left\vert \log\frac{\widetilde{t}_{u}}{u}+\log\left(
a(x)+\frac{(1-\widetilde{t}_{u})u}{\widetilde{t}_{u}}\right)  \right\vert \\
&  <\max(\left\vert \log\frac{H}{2}\right\vert ,\text{ }\left\vert \log
\frac{3H}{2}\right\vert )+\log2+\left\vert \log a(x)\right\vert \text{.}%
\end{align*}

(2) If $a(x)<1/2$,
\[
\left\vert \log\left(  a(x)\widetilde{t}_{u}/u-\widetilde{t}_{u}+1\right)
\right\vert <\max(\left\vert \log\frac{H}{2}\right\vert ,\text{ }\left\vert
\log\frac{3H}{2}\right\vert )+\left\vert \log a(x)\right\vert \text{.}%
\]
Using the above properties, we can apply Lebesgue theorem as follows:%
\begin{align*}
\lim_{u\rightarrow(1/H)^{+}}\widetilde{G}_{u}(\widetilde{t}_{u})  &
=\lim_{u\rightarrow(1/H)^{+}}\exp\left(  \int_{I}\log\left(  a(x)\widetilde
{t}_{u}/u-\widetilde{t}_{u}+1\right)  d(F(x))\right) \\
&  =H\exp\left(  \int_{I}\log a(x)d(F(x))\right)  .
\end{align*}

\hfill$\square$

\bigskip

\noindent\textbf{Lemma 4.17}. \textit{If }$\xi=0$\textit{, }$H<+\infty
,$\textit{ and }$\widetilde{G}<+\infty$\textit{, }$\lim_{u\rightarrow
(1/H)^{-}}\widetilde{G}_{u}(1^{-})$ $=H\exp\left(  \int_{I}\log
a(x)d(F(x))\right)  .$

\bigskip

\noindent\textbf{Proof.} In the case in which $\xi=0$, based on the definition
which is mentioned beneath the proof of Lemma 4.10, we have $\widetilde{G}%
_{u}(1^{-})=\exp\left(  \int_{I}\log a(x)d(F(x))\right)  /u$. Thus, we obtain
the conclusion.\hfill$\square$

\bigskip

\noindent\textbf{Lemma 4.18}.\textit{ If }$\xi=0$\textit{, }$H<+\infty
,$\textit{ and }$\widetilde{G}<+\infty$\textit{, }$\lim_{u\rightarrow0^{+}%
}\widetilde{G}_{u}(1^{-})$ $=+\infty.$

\bigskip

\textbf{Proof.} Due to Lemmas 4.11, 4.12, and 4.16, we have $H$ $\exp(\int
_{I}\log a(x)d(F(x)))$ $>1$. Therefore, by the definition of $\widetilde
{G}_{u}(1^{-})$\ we obtain%
\[
\lim_{u\rightarrow0^{+}}\widetilde{G}_{u}(1^{-})=\lim_{u\rightarrow0^{+}}%
\exp\left(  \int_{I}\log a(x)d(F(x))\right)  /u\geq\lim_{u\rightarrow0^{+}%
}1/H/u=+\infty.
\]

\hfill$\square$

\bigskip

\noindent\textbf{Lemma 4.19.} \textit{If }$\xi>0$\textit{ and }$H_{\xi
}=+\infty$\textit{, }$\lim_{u\rightarrow\xi^{+}}\widetilde{G}_{u}%
(\widetilde{t}_{u})=+\infty.$

\bigskip

\noindent\textbf{Proof.} Lemma 4.12 ensures the existence of $\lim
_{u\rightarrow\xi^{+}}\widetilde{G}_{u}(\widetilde{t}_{u})$, which is finite
or $+\infty$. If $a(x)>\xi$, then $(a(x)-\xi)/(2(u-\xi))$ is strictly
decreasing with respect to $u\in(\xi$, $E).$ Using Lebesgue theorem, we have
\begin{align*}
\lim_{u\rightarrow\xi^{+}}\widetilde{G}_{u}(\widetilde{t}_{u})  &
\geq\underline{\lim}_{u\rightarrow\xi^{+}}\widetilde{G}_{u}(\frac{u}{2(u-\xi
)})\\
&  \geq\underline{\lim}_{u\rightarrow\xi^{+}}\frac{1}{2}\exp\left(
\int_{a(x)>\xi}\log\left(  \frac{a(x)-\xi}{2(u-\xi)}+\frac{1}{2}\right)
d(F(x))\right)  =+\infty.
\end{align*}

\hfill$\square$

\bigskip

\noindent\textbf{Lemma 4.20.} \textit{If }$\xi>0$\textit{, }$H_{\xi}<+\infty
,$\textit{ and }$\widetilde{G}<+\infty$\textit{, }$\lim_{u\rightarrow
(\xi+1/H_{\xi})^{+}}\widetilde{G}_{u}(\widetilde{t}_{u})$ $=H_{\xi}\exp
(\int_{I}\log\left(  a(x)-\xi\right)  d(F(x))).$

\bigskip

\noindent\textbf{Proof.} An argument similar to that in the proof of Lemma
4.16 ensures that $\log(a(x)-\xi)$\ is integrable.

It should be noted that $\lim_{u\rightarrow(\xi+1/H_{\xi})^{+}}\widetilde
{t}_{u}=1+\xi H_{\xi}$ (Lemma 3.9). From the fact that $\lim_{u\rightarrow
(\xi+1/H_{\xi})^{+}}\widetilde{t}_{u}/u=H_{\xi}$ and $\lim_{u\rightarrow
(\xi+1/H_{\xi})^{+}}(\xi\widetilde{t}_{u}-u\widetilde{t}_{u}+u)/\widetilde
{t}_{u}=0,$ we can choose $0<\delta<\min(\xi+1/H_{\xi},$ $E-\xi-1/H_{\xi})$,
such that $H_{\xi}/2<\widetilde{t}_{u}/u<3H_{\xi}/2$ and $(\xi\widetilde
{t}_{u}-u\widetilde{t}_{u}+u)/\widetilde{t}_{u}<1/2$ for each $u\in
(\xi+1/H_{\xi},$ $\xi+1/H_{\xi}+\delta)$. Therefore, we have the following properties.

(1) If $a(x)\geq\xi+1/2$,
\begin{align*}
&  \left\vert \log\left(  a(x)\widetilde{t}_{u}/u-\widetilde{t}_{u}+1\right)
\right\vert \\
&  =\left\vert \log\frac{\widetilde{t}_{u}}{u}+\log\left(  a(x)-\xi+\frac
{\xi\widetilde{t}_{u}-u\widetilde{t}_{u}+u}{\widetilde{t}_{u}}\right)
\right\vert \\
&  <\max(\left\vert \log\frac{H_{\xi}}{2}\right\vert ,\text{ }\left\vert
\log\frac{3H_{\xi}}{2}\right\vert )+\log2+\left\vert \log\left(
a(x)-\xi\right)  \right\vert \text{.}%
\end{align*}

(2) If $a(x)<\xi+1/2$,
\[
\left\vert \log\left(  a(x)\widetilde{t}_{u}/u-\widetilde{t}_{u}+1\right)
\right\vert <\max(\left\vert \log\frac{H_{\xi}}{2}\right\vert ,\text{
}\left\vert \log\frac{3H_{\xi}}{2}\right\vert )+\left\vert \log\left(
a(x)-\xi\right)  \right\vert \text{.}%
\]
Using the above properties, we can apply Lebesgue theorem as follows.%
\begin{align*}
\lim_{u\rightarrow(\xi+1/H_{\xi})^{+}}\widetilde{G}_{u}(\widetilde{t}_{u})  &
=\lim_{u\rightarrow(\xi+1/H_{\xi})^{+}}\exp\left(  \int_{I}\log\left(
a(x)\widetilde{t}_{u}/u-\widetilde{t}_{u}+1\right)  d(F(x))\right) \\
&  =H_{\xi}\exp\left(  \int_{I}\log(a(x)-\xi)d(F(x))\right)  .
\end{align*}

\hfill$\square$

\bigskip

\noindent\textbf{Lemma 4.21.}\textit{ If }$\xi>0$\textit{, }$\widetilde
{G}_{1/H}(\widetilde{t}_{1/H})=H\exp\left(  \int_{I}\log a(x)d(F(x))\right)
.$

\bigskip

\noindent\textbf{Proof.} It should be noted that $0<\xi<1/H<E$ and
$1<1/H/(1/H-\xi)$. From Lemma 3.8, we have $\widetilde{t}_{1/H}=1$. Thus,
\begin{align*}
\widetilde{G}_{1/H}(\widetilde{t}_{1/H})  &  =\exp\left(  \int_{I}\log\left(
a(x)/(1/H)-1+1\right)  d(F(x))\right) \\
&  =H\exp\left(  \int_{I}\log a(x)d(F(x))\right)  .
\end{align*}

\hfill$\square$

\bigskip

\noindent\textbf{Lemma 4.22.} \textit{If }$\xi>0$\textit{, }$\lim
_{u\rightarrow(1/H)^{-}}$\textit{ }$\widetilde{G}_{u}(1)=H\exp\left(  \int
_{I}\log a(x)d(F(x))\right)  .$

\bigskip

\noindent\textbf{Proof. }From $0<\xi<u$, we have $1<u/(u-\xi)$. Thus, by Lemma
4.21 we obtain
\[
\lim_{u\rightarrow(1/H)^{-}}\widetilde{G}_{u}(1)=\widetilde{G}_{1/H}%
(1)=H\exp\left(  \int_{I}\log a(x)d(F(x))\right)  .
\]

\hfill$\square$

\bigskip

If $\xi>0$ and $0<u\leq\xi$, $a(x)t/u-t+1\geq1$ for each $t>0$. Therefore, we
can expand the definition of $\widetilde{G}_{u}(t)=\exp(\int_{I}\log\left(
a(x)t/u-t+1\right)  d(F(x)))$, which is greater than $1$ and finite or
$+\infty$, in the domain $0<u\leq\xi$ and $t>0$.

\bigskip

\noindent\textbf{Lemma 4.23.} \textit{If }$\xi>0$\textit{, }$\lim
_{u\rightarrow0^{+}}$\textit{ }$\widetilde{G}_{u}(1)=+\infty.$

\bigskip

\noindent\textbf{Proof. }It should be noted that $\widetilde{G}_{1/H}%
(\widetilde{t}_{1/H})=H\exp\left(  \int_{I}\log a(x)d(F(x))\right)  >1$ (Lemma
4.21). From the expansion of $\widetilde{G}_{u}(t)$ which is defined beneath
the proof of Lemma 4.22, we have%
\[
\lim_{u\rightarrow0^{+}}\widetilde{G}_{u}(1)=\lim_{u\rightarrow0^{+}}%
\exp\left(  \int_{I}\log a(x)d(F(x))\right)  /u\geq\lim_{u\rightarrow0^{+}%
}1/H/u=+\infty.
\]

\hfill$\square$

\bigskip

\noindent\textbf{Lemma 4.24.} \textit{If }$\xi<0$\textit{, }$\xi+1/H_{\xi}%
>0$\textit{ and }$\widetilde{G}<+\infty$\textit{, }$\lim_{u\rightarrow
(\xi+1/H_{\xi})^{+}}\widetilde{G}_{u}(\widetilde{t}_{u})$ $=H_{\xi}\exp\left(
\int_{I}\log\left(  a(x)-\xi\right)  d(F(x))\right)  .$

\bigskip

\noindent\textbf{Proof}. It should be noted that $H_{\xi}<+\infty$ and
$\lim_{u\rightarrow(\xi+1/H_{\xi})^{+}}\widetilde{t}_{u}=1+\xi H_{\xi}%
$\ (Lemma 3.10). The proof is formally the same as that of Lemma
4.20.\hfill$\square$

\bigskip

\noindent\textbf{Lemma 4.25.} \textit{If }$\xi<0$\textit{, }$\xi+1/H_{\xi}%
>0,$\textit{ and }$\widetilde{G}<+\infty$\textit{, }$\lim_{u\rightarrow
(\xi+1/H_{\xi})^{-}}$

$\widetilde{G}_{u}((u/(u-\xi))^{-})=H_{\xi}\exp\left(  \int_{I}\log\left(
a(x)-\xi\right)  d(F(x))\right)  .$

\bigskip

\noindent\textbf{Proof}. We obtain the conclusion using the definition which
is mentioned beneath the proof of Lemma 4.10.\hfill$\square$

\bigskip

\noindent\textbf{Lemma 4.26.} \textit{If }$\xi<0$\textit{, }$\xi+1/H_{\xi}%
>0,$\textit{ and }$\widetilde{G}<+\infty$\textit{, }$\lim_{u\rightarrow0^{+}%
}\widetilde{G}_{u}(u/(u-\xi)^{-})$ $=\exp\left(  \int_{I}\log\left(
a(x)-\xi\right)  d(F(x))\right)  $ $/(-\xi).$

\bigskip

\noindent\textbf{Proof}. We obtain the conclusion by applying the same process
as in Lemma 4.25.

\hfill$\square$

\bigskip

From Lemma 3.15, $\widetilde{t}_{u}/u$ is strictly decreasing with respect to
$u\in(0$, $E)$, if $\xi+1/H_{\xi}\leq0$. Therefore, $\lim_{u\rightarrow0^{+}%
}\widetilde{t}_{u}/u$ exists, and we denote it by $\eta>0$. From
$0<\widetilde{t}_{u}<u/(u-\xi)$, we have $\eta\leq-1/\xi$.

\bigskip

\noindent\textbf{Lemma 4.27.} \textit{If }$\xi<0$\textit{ and }$\xi+1/H_{\xi
}<0$\textit{, }$\eta<-1/\xi$\textit{.}

\bigskip

\noindent\textbf{Proof}. It should be noted that the definition $w_{u}%
(\widetilde{t}_{u})=0$ implies that%
\[
\int_{I}\frac{1}{a(x)\frac{\widetilde{t}_{u}}{u}-\widetilde{t}_{u}%
+1}d(F(x))=1.
\]
From Lemma 3.12, we have $\lim_{u\rightarrow0^{+}}\widetilde{t}_{u}=0$. Using
Fatou's lemma, we obtain%
\[
\int_{I}\frac{1}{a(x)\eta+1}d(F(x))=\int_{I}\underline{\lim}_{u\rightarrow
0^{+}}\frac{1}{a(x)\frac{\widetilde{t}_{u}}{u}-\widetilde{t}_{u}+1}%
d(F(x))\leq1.
\]
Assume that $\eta=-1/\xi$, we have%
\[
\int_{I}\frac{1}{-a(x)/\xi+1}d(F(x))=-\xi H_{\xi}\leq1.
\]
This implies that $\xi+1/H_{\xi}\geq0$, which is a contradiction.\hfill
$\square$

\bigskip

\noindent\textbf{Lemma 4.28.} \textit{If }$\xi<0$\textit{ and }$\xi+1/H_{\xi
}=0$\textit{, }$\eta=-1/\xi$.

\bigskip

\noindent\textbf{Proof.} Since $H_{\xi}<+\infty$, $1/(a(x)-\xi)$ is
integrable. Thus, from $\int_{I}1/(a(x)\widetilde{t}_{u}/u-\widetilde{t}_{u}$
$+1)d(F(x))=1$ and by Lebesgue theorem, we have $\int_{I}1/(a(x)\eta
+1)d(F(x))=1$. On the other hand, it is clear that $\int_{I}1/(a(x)\times
0+1)d(F(x))=1$ and $\int_{I}1/(a(x)(-1/\xi)$ $+1)d(F(x))=-\xi H_{\xi}=1$. This
implies that the equation $\int_{I}1/(a(x)y+1)d(F(x))$ $=0$ with respect to
$y\in\lbrack0,$ $-1/\xi]$ has three solutions $y=0,$ $y=\eta,$ and $y=-1/\xi$.
Note that%
\[
\frac{\partial^{2}}{\partial y^{2}}\int_{I}\frac{1}{a(x)y+1}d(F(x))=\int
_{I}\frac{a(x)^{2}}{\left(  a(x)y+1\right)  ^{3}}d(F(x))>0.
\]
Therefore, the equation $\int_{I}1/(a(x)y+1)d(F(x))=0$\ has at most two
solutions. This implies that $\eta=-1/\xi$.\hfill$\square$

\bigskip

\noindent\textbf{Lemma 4.29.} \textit{If }$\xi<0$\textit{, }$\xi+1/H_{\xi}%
\leq0,$\textit{ and }$\widetilde{G}<+\infty$\textit{, }$\lim_{u\rightarrow
0^{+}}\widetilde{G}_{u}(\widetilde{t}_{u})$ $=\exp\left(  \int_{I}%
\log(a(x)\eta+1)d(F(x))\right)  .$

\bigskip

\noindent\textbf{Proof}. Lemma 3.12 implies that $\lim_{u\rightarrow0^{+}%
}\widetilde{t}_{u}=0$. From Lemma 3.15, $\widetilde{t}_{u}/u$ is strictly
decreasing with respect to $u\in(0$, $E)$. Due to Lemma 3.16, $\widetilde
{t}_{u}$ is strictly increasing with respect to $u\in(0$, $u_{\max})$.
Therefore, if $a(x)>0$ then $a(x)\frac{\widetilde{t}_{u}}{u}-\widetilde{t}%
_{u}+1$ is strictly decreasing with respect to $u\in(0$, $u_{\max})$. This
ensures that%
\[
\lim_{u\rightarrow0^{+}}\int_{a(x)>0}\log(\frac{a(x)}{u}\widetilde{t}%
_{u}-\widetilde{t}_{u}+1)d(F(x))=\int_{a(x)>0}\log(a(x)\eta+1)d(F(x)).
\]

If $H_{\xi}<+\infty,$ using Jensen's inequality theorem, we see that
$\log(a(x)-\xi)$ is integrable. If $a(x)\leq0$ and $0<u<\min(-\xi$, $E)$, then
we have%
\[
0\leq\frac{a(x)-\xi}{-2\xi}\leq\frac{a(x)-\xi}{u-\xi}<a(x)\frac{\widetilde
{t}_{u}}{u}-\widetilde{t}_{u}+1<1,
\]
and%
\[
\left\vert \log\left(  a(x)\widetilde{t}_{u}/u-\widetilde{t}_{u}+1\right)
\right\vert <\left\vert \log\left(  a(x)-\xi\right)  \right\vert +\left\vert
\log(-2\xi)\right\vert \text{.}%
\]
Therefore, we can apply Lebesgue theorem to the following equality.%
\[
\lim_{u\rightarrow0^{+}}\int_{a(x)\leq0}\log(\frac{a(x)}{u}\widetilde{t}%
_{u}-\widetilde{t}_{u}+1)d(F(x))=\int_{a(x)\leq0}\log(a(x)\eta+1)d(F(x)).
\]
Thus, we accomplish
\[
\lim_{u\rightarrow0^{+}}\widetilde{G}_{u}(\widetilde{t}_{u})=\exp\left(
\int_{I}\log(a(x)\eta+1)d(F(x))\right)  .
\]

If $H_{\xi}=+\infty,$ from Lemma 4.27, $\eta<-1/\xi$. If we assign
$\varepsilon:=(\xi\eta+1)$ $/2>0$, then, there exists $\delta>0$ such that
$\widetilde{t}_{u}<\varepsilon$ for each $u\in(0$, $\delta)$. Hence, if
$a(x)\leq0$, then we have%
\[
1>\frac{a(x)}{u}\widetilde{t}_{u}-\widetilde{t}_{u}+1\geq\xi\eta
-\varepsilon+1=\varepsilon.
\]
Thus, we can apply Lebesgue theorem in the domain $\left\{  \text{ }x\text{
}|\text{ }a(x)<0\text{ }\right\}  $ and obtain the conclusion.\hfill$\square$

\bigskip

Here, we redefine $\eta$ to be $\lim_{u\rightarrow0^{+}}\overline{t}_{u}/u$.
If $\xi<0$ and $\xi+1/H_{\xi}>0$, by the definition of $\overline{t}_{u}$, we
have $\eta=-1/\xi$. Therefore, Lemma 4.26 implies that $\lim_{u\rightarrow
0^{+}}\widetilde{G}_{u}(\overline{t}_{u})=\exp\left(  \int_{I}\log\left(
a(x)\eta+1\right)  d(F(x))\right)  $. Summing up the above-mentioned Lemmas,
we obtain the following

\bigskip

\noindent\textbf{Theorem 4.1.} \textit{If }$\widetilde{G}<+\infty$\textit{,
}$\widetilde{G}_{u}(\overline{t}_{u})$\textit{ is continuous and strictly
decreasing with respect to }$u\in(0$, $E)$\textit{. The range of }%
$\widetilde{G}_{u}(\overline{t}_{u})$\textit{ is }$(1,$\textit{ }$+\infty
)$\textit{ }(\textit{if }$\xi\geq0$)\textit{ or }$(1,$\textit{ }$\exp\left(
\int_{I}\log\left(  a(x)\eta+1\right)  d(F(x))\right)  )$\textit{ }(\textit{if
}$\xi<0$).

\bigskip

Afterward, we will show the equality $\widetilde{G}_{u}(\overline{t}%
_{u})=G_{u}(t_{u})$ (Theorem 5.1).

\bigskip

\noindent{\large \textbf{5. Double sequence of random variables}}

It should be noted that a series of step functions exists such that $\xi\leq
f_{N}(x)$ $\leq f_{N+1}(x)\leq a(x)$ and $\lim_{N\rightarrow+\infty}%
f_{N}(x)=a(x)$ for each $x\in I,$ in which $\xi$ $=\inf_{x\in I}%
\,a(x)>-\infty$ is the essential infimum.

For example, for each positive integer $N$, assign $M:=2^{N}N+1$ and
\begin{equation}
f_{N}(x):=\left\{
\begin{array}
[c]{cc}%
a_{j}:=\xi+\frac{j-1}{2^{N}}, & \text{if }\xi+\frac{j-1}{2^{N}}\leq
a(x)<\xi+\frac{j}{2^{N}}\text{ \ }(1\leq j\leq M-1),\\
a_{M}:=\xi+N, & \text{if }a(x)\geq\xi+N.
\end{array}
\right.  \tag{12}%
\end{equation}

In general, suppose $\{f_{N}(x)$ $|$ $x\in I\}=\{a_{j}$ $|$ $j=1,...,M\}$,
where $\xi=a_{1}$ $<a_{2}<\cdots<a_{M}<+\infty$. Set $p_{j}:=\int_{a_{j}\leq
a(x)<a_{j+1}}d(F(x))$ \ $(1\leq j\leq M-1)$ and $p_{M}:=\int_{a(x)\geq a_{M}%
}d(F(x))$, then we have $\sum_{j=1}^{M}p_{j}=1$.

Assume $u>0$, $t\geq0,$ and $\xi t/u-t+1>0$. Then $a_{j}t/u-t+1>0$ for each
$1\leq j\leq M$. For the game $\{(a_{j},$ $p_{j})$ $|$ $j=1,2,\cdots,M\}$, the
growth rate per attempt after $n$ attempts is
\begin{equation}
\left(  \prod\limits_{k=1}^{s}\left(  a_{j_{k}}t/u-t+1\right)  ^{m_{j_{k}}%
}\right)  ^{\frac{1}{n}}, \tag{13}%
\end{equation}
where $a_{j_{k}}$ occurs $m_{j_{k}}$ times $(m_{j_{1}}+m_{j_{2}}%
+\cdots+m_{j_{s}}=n)$ with probability $p_{j_{1}}^{m_{j_{1}}}\cdots p_{j_{s}%
}^{m_{j_{s}}}$. Such event has $n!$ $/(m_{j_{1}}!m_{j_{2}}!\cdots m_{j_{s}}%
!)$\ permutation patterns. We denote $X_{N,n}$ by this random variable. Then,
the expectation $E[X_{N,\text{ }n}]$ is expressed as%
\begin{align}
&  \sum\limits_{m_{j_{1}}+\cdots+m_{j_{s}}=n}\frac{n!}{m_{j_{1}}!\cdots
m_{j_{s}}!}\left(  \prod\limits_{k=1}^{s}\left(  a_{j_{k}}t/u-t+1\right)
^{m_{j_{k}}}\right)  ^{\frac{1}{n}}p_{j_{1}}^{m_{j_{1}}}\cdots p_{j_{s}%
}^{m_{j_{s}}}\tag{14}\\
&  =\left(  \sum\limits_{j=1}^{M}\left(  a_{j}t/u-t+1\right)  ^{\frac{1}{n}%
}p_{j}\right)  ^{n}\nonumber\\
&  =\exp\left(  \frac{\log\left(  \sum\limits_{j=1}^{M}\left(  a_{j}%
t/u-t+1\right)  ^{\frac{1}{n}}p_{j}\right)  }{\frac{1}{n}}\right)  .\nonumber
\end{align}
Moreover, the variance $V[X_{N,n}]$ is expressed as%
\begin{align}
&  \sum\limits_{m_{j_{1}}+\cdots+m_{j_{s}}=n}\frac{n!}{m_{j_{1}}!\cdots
m_{j_{s}}!}\left(  \prod\limits_{k=1}^{s}\left(  a_{j_{k}}t/u-t+1\right)
^{m_{j_{k}}}\right)  ^{\frac{2}{n}}p_{j_{1}}^{m_{j_{1}}}\cdots p_{j_{s}%
}^{m_{j_{s}}}\tag{15}\\
&  -E[X_{N,n}]^{2}\nonumber\\
&  =\exp\left(  \frac{\log\left(  \sum\limits_{j=1}^{M}\left(  a_{j}%
t/u-t+1\right)  ^{\frac{2}{n}}p_{j}\right)  }{\frac{1}{n}}\right) \nonumber\\
&  -\exp\left(  \frac{2\log\left(  \sum\limits_{j=1}^{M}\left(  a_{j}%
t/u-t+1\right)  ^{\frac{1}{n}}p_{j}\right)  }{\frac{1}{n}}\right)  .\nonumber
\end{align}

\bigskip

\noindent\textbf{Lemma 5.1.}%
\[
\lim_{N\rightarrow+\infty}\left(  \lim_{n\rightarrow+\infty}E[X_{N,n}]\right)
=\exp\left(  \int_{I}\log\left(  a(x)t/u-t+1\right)  d(F(x))\right)  .
\]

\bigskip

\noindent\textbf{Proof.} If $n$ approaches $+\infty$, using L'Hopital's
theorem, we obtain%
\begin{align*}
\lim_{n\rightarrow+\infty}E[X_{N,n}]  &  =\lim_{h\rightarrow0^{+}}\exp\left(
\frac{\sum\limits_{j=1}^{M}\left(  a_{j}t/u-t+1\right)  ^{h}\log\left(
a_{j}t/u-t+1\right)  p_{j}}{\sum\limits_{j=1}^{M}\left(  a_{j}t/u-t+1\right)
^{h}p_{j}}\right) \\
&  =\exp\left(  \int_{I}\log\left(  f_{N}(x)t/u-t+1\right)  d(F(x))\right)  .
\end{align*}
Since $f_{N}(x)\leq f_{N+1}(x)\leq a(x),$ and $\lim_{N\rightarrow+\infty}%
f_{N}(x)=a(x)$, we obtain the conclusion.

\hfill$\square$

\bigskip

It is easily verified that if $\widetilde{G}<+\infty$, $\lim_{n\rightarrow
+\infty}V[X_{N,n}]=0$ for each $N$.

\bigskip

\noindent\textbf{Lemma 5.2.} \textit{If }$\nu>0$\textit{, the following
statements are equivalent.}

(1) $\int_{a(x)>1}a(x)^{\nu}d(F(x))<+\infty$.

(2) $\int_{I}(a(x)t/u-t+1)^{^{\nu}}d(F(x))<+\infty$ \textit{for each }%
$u>0$\textit{, }$t>0,$\textit{ and} $\xi t/u-t+1>0.$

(3) $\int_{I}(a(x)t_{1}/u_{1}-t_{1}+1)^{^{\nu}}d(F(x))<+\infty$ \textit{for
some }$u_{1}>0$\textit{, }$t_{1}>0,$\textit{ and} $\xi t_{1}/u_{1}-t_{1}+1>0.$

\bigskip

\noindent\textbf{Proof.} If $a(x)>2u/t\times\left\vert 1-t\right\vert $, we
have $a(x)<2u/t\times(a(x)t/u-t+1)$ and $a(x)t/u-t+1<2t/u\times a(x).$ This
implies the conclusion.\hfill$\square$

\bigskip

We say that a game $(a(x)$, $F(x))$ is \textit{effective} when $\int
_{a(x)>1}a(x)^{\nu}d(F(x))<+\infty$ for some $\nu>0$, with the additional
conditions $E>0$ and $\xi>-\infty$.

It should noted that for each $0<\nu<1,$ there exists $h_{\nu}$ such that
$\log(x+1)$ $<h_{\nu}x^{\nu}$ for each $x>0$.

\bigskip

\noindent\textbf{Lemma 5.3.} If a game is effective,\ $\widetilde{G}<+\infty$.

\bigskip

\noindent\textbf{Proof.} If $a(x)>1$ and $0<\nu<1$, we have $\log a(x)<h_{\nu
}(a(x)-1)^{\nu}<h_{\nu}a(x)^{\nu}$. If $a(x)>1$ and $\nu\geq1$, we have $\log
a(x)<a(x)\leq a(x)^{\nu}$. Thus, we obtain the conclusion.\hfill$\square$

\bigskip

For example, the game $(\exp(1/\sqrt{x}),$ $x)$ $(x\in(0$, $1])$ is
ineffective and $\widetilde{G}<+\infty$.

\bigskip

\noindent\textbf{Lemma 5.4.} \textit{If a game is effective,\ }$\lim
_{N\rightarrow+\infty}\left(  \lim_{n\rightarrow+\infty}E[X_{N,n}]\right)  $

$=\lim_{n\rightarrow+\infty}\left(  \lim_{N\rightarrow+\infty}E[X_{N,n}%
]\right)  =\lim_{\substack{n\rightarrow+\infty\\N\rightarrow+\infty}%
}E[X_{N,n}]<+\infty.$

\textit{If a game is ineffective, }$\lim_{N\rightarrow+\infty}E[X_{N,n}%
]=+\infty$ \textit{for each} $n$, \textit{if} $t>0.$

\bigskip

\noindent\textbf{Proof.} For the assumption $\int_{a(x)>1}a(x)^{\nu
}d(F(x))<+\infty$, we can assume that $\nu<1$. Suppose $0<h<\nu/2$ and
$a(x)>u$, we have%
\begin{align*}
\left\vert \frac{\partial}{\partial h}(a(x)t/u-t+1)^{h}\right\vert  &
<(a(x)t/u-t+1)^{\nu/2}\log(a(x)t/u-t+1)\\
&  <h_{\nu/2}(a(x)t/u-t+1)^{\nu}.
\end{align*}
This guarantees that
\[
\frac{\partial}{\partial h}\int_{I}(a(x)t/u-t+1)^{h}d(F(x))=\int_{I}%
\frac{\partial}{\partial h}(a(x)t/u-t+1)^{h}d(F(x))<+\infty.
\]
Thus, since $f_{N}(x)\leq f_{N+1}(x)\leq a(x),$ and $\lim_{N\rightarrow
+\infty}f_{N}(x)=a(x)$, we have%
\begin{align*}
&  \lim_{n\rightarrow+\infty}\left(  \lim_{N\rightarrow+\infty}E[X_{N,n}%
]\right) \\
&  =\exp\left(  \lim_{h\rightarrow0^{+}}\frac{\int_{I}(a(x)t/u-t+1)^{h}%
\log(a(x)t/u-t+1)d(F(x))}{\int_{I}(a(x)t/u-t+1)^{h}d(F(x))}\right) \\
&  =\exp\left(  \int_{I}\log\left(  a(x)t/u-t+1\right)  d(F(x))\right)
<+\infty.
\end{align*}
Using Lemma 5.1, we have
\[
\lim_{N\rightarrow+\infty}\left(  \lim_{n\rightarrow+\infty}E[X_{N,n}]\right)
=\lim_{n\rightarrow+\infty}\left(  \lim_{N\rightarrow+\infty}E[X_{N,n}%
]\right)  .
\]
Set $\alpha:=\lim_{N\rightarrow+\infty}\left(  \lim_{n\rightarrow+\infty
}E[X_{N,n}]\right)  =\lim_{n\rightarrow+\infty}\left(  \lim_{N\rightarrow
+\infty}E[X_{N,n}]\right)  $, $U_{N}$ $:=\lim_{n\rightarrow+\infty}%
E[X_{N,n}],$ and $W_{n}:=\lim_{N\rightarrow+\infty}E[X_{N,n}]$. Since
$f_{N}(x)\leq f_{N+1}(x)\leq a(x),$ $E[X_{N,n}]$ increases with respect to
$N.$ By setting $h:=1/n$ and applying L'Hopital's theorem twice, we have
\begin{align*}
&  \lim_{h\rightarrow0^{+}}\frac{\partial}{\partial h}E[X_{N,n}]\\
&  =\frac{\exp\left(  \int_{I}\log(f_{N}(x)t/u-t+1)d(F(x))\right)  }{2}%
\times\\
&  \left(
\begin{array}
[c]{c}%
\int_{I}\left(  \log(f_{N}(x)t/u-t+1)\right)  ^{2}d(F(x))\\
-\left(  \int_{I}\log(f_{N}(x)t/u-t+1)d(F(x))\right)  ^{2}%
\end{array}
\right)  \geq0.
\end{align*}
This implies that $E[X_{N,n}]$ decreases with sufficiently large $n$. For
$\varepsilon>0$, there exists $N_{0}>0$ such that $\left\vert U_{N}%
-\alpha\right\vert <\varepsilon$ for each $N\geq N_{0}.$ Moreover, for
$\varepsilon>0,$ there exists $N_{1}>0$ such that $\left\vert W_{n}%
-\alpha\right\vert <\varepsilon$ for each $n\geq N_{1}.$ Therefore, we have
$\alpha-\varepsilon<U_{N}\leq E[X_{N,n}]\leq W_{n}<\alpha+\varepsilon,$ which
implies $\lim_{\substack{n\rightarrow+\infty\\N\rightarrow+\infty}%
}E[X_{N,n}]=\alpha$.

If $\int_{a(x)>1}a(x)^{\nu}d(F(x))=+\infty$ for each $\nu>0$, using Lemma 5.2,
we have%
\[
\lim_{N\rightarrow+\infty}E[X_{N,n}]=\exp\left(  \frac{\log\left(  \int
_{I}(a(x)t/u-t+1)^{\frac{1}{n}}d(F(x))\right)  }{\frac{1}{n}}\right)  =+\infty
\]
for each $n\geq1$ and $t>0$.\hfill$\square$

\bigskip

\noindent\textbf{Lemma 5.5.} \textit{If a game is effective,\ }$\lim
_{\substack{n\rightarrow+\infty\\N\rightarrow+\infty}}V[X_{N,n}]=0.$

\bigskip

\noindent\textbf{Proof.} This can be proved in a manner similar to Lemma
5.4.\hfill$\square$

\bigskip

\noindent\textbf{Lemma 5.6.} \textit{If a game is effective,\ }%
\[
\lim_{\substack{n\rightarrow+\infty\\N\rightarrow+\infty}}E[\left(
X_{N,n}-\exp\left(  \int_{I}\log\left(  a(x)t/u-t+1\right)  d(F(x))\right)
\right)  ^{2}]=0.
\]

\bigskip

\noindent\textbf{Proof.} The equality $E[\left(  X_{N,n}-c\right)
^{2}]=V[X_{N,n}]+\left(  E[X_{N,n}]-c\right)  ^{2}$ for each $c$ implies the
conclusion (Lemmas 5.1, 5.4, and 5.5).\hfill$\square$

\bigskip

We denote $\exp(\int_{I}\log\left(  a(x)t/u-t+1\right)  d(F(x)))$ by
$G_{u}(t)$ for each $u>0$, $0$ $\leq t\leq1$, and $\xi t/u-t+1>0$. We refer to
this as \textit{the limit expectation of growth rate}. We adopt the criteria
$\sup_{0\leq t\leq1,\text{ }\xi t-t+1>0}$ $G_{1}(t)$. In order to explain an
advantage of this criteria, we compare the following two games. After the
consideration below, it is reasonable for investors to prefer Game-2, which
has the lower expectation $1.3125$ $(<1.5)$ than that of Game-1.

\bigskip

\noindent\textbf{Game-1}. The profit $3$ or $0$ occurs with probability $0.5$.
In this case, we have $E=1.5$, $\max_{0\leq t\leq1,\text{ }\xi t-t+1>0}$
$G_{1}(t)=\sqrt{9/8}\fallingdotseq1.0607,$ and $\widetilde{t}_{1}=0.25$
(Example 6.3).

If the investors continue to invest all their current capital in Game-1 with
price $1$, then, after $30$ attempts they will be ruined with probability
$1-0.5^{30}$ $\fallingdotseq0.99999999907$.

If the investors continue to invest $0.25$ of their current capital in Game-1
with price $1$, after $30$ attempts their capital will increase by a factor of
$E[X_{N,\text{ }30}]$ $=((3\times0.25+0.75)^{1/30}/2+0.75^{1/30}%
/2)^{30}\fallingdotseq1.0628$ $(1.0628^{30}\fallingdotseq6.2165)$ with a small
variance $V[X_{N,\text{ }30}]\fallingdotseq0.0045$ for each $N\geq3$, where
$f_{N}(x)$\ is defined by (12). Moreover, $E[X_{N,n}]<1.0929$ for each
$n\geq2$ and $N\geq3.$

\bigskip

\noindent\textbf{Game-2.} The profit $2$ or $0.625$ occurs with probability
$0.5$. In this case, we have $E$ $\fallingdotseq1.3125$, $\max_{0\leq
t\leq1,\text{ }\xi t-t+1>0}$ $G_{1}(t)=\sqrt{121/96}\fallingdotseq1.1227,$ and
$\widetilde{t}_{1}=5/6\fallingdotseq0.8333$ (Example 6.3).

If the investors continue to invest $0.8333$ of their current capital in
Game-2 with price $1$, after $30$ attempts their capital will increase by a
factor of $E[X_{N,\text{ }30}]$ $=((2\times0.8333+0.1667)^{1/30}%
/2+(0.625\times0.8333+0.1667)^{1/30}/2)^{30}$ $\fallingdotseq1.1272$
$(1.1272^{30}$ $\fallingdotseq36.3103)$ with a small variance $V[X_{N,\text{
}30}]\fallingdotseq0.0102$ for each $N\geq3$, where $f_{N}(x)$\ is defined by
(12). Moreover, $E[X_{N,n}]>1.1227$ for each $n\geq2$ and $N\geq3.$

\bigskip

In general, the asymptotic optimality principle states that any alternative is
dominated in the long run by the log-optimum strategy (Algoet and Cover
(1988)). However, it should be noted that, despite Lemma 5.1, when a game is
ineffective, the limit expectation of growth rate dose not have a solid
significance as shown by Lemma 5.4.

\bigskip

We say that $t_{u}$ is \textit{the optimal proportion of investment} in order
to continue to invest without borrowing with respect to $u>0,$ if
\begin{equation}
\overline{\lim}_{\substack{\rho\rightarrow t_{u}\\0\leq\rho\leq1\\\xi
\rho/u-\rho+1>0}}\int_{I}\log\frac{a(x)t/u-t+1}{a(x)\rho/u-\rho+1}d(F(x))\leq0
\tag{16}%
\end{equation}
for each $0\leq t\leq1$ and $\xi t/u-t+1>0.$ It follows that $0\leq t_{u}%
\leq1$ and $\xi t_{u}/u-t_{u}+1$ $\geq0$.

It is clear that $G_{u}(t)=\widetilde{G}_{u}(t)$ on the intersection of
domains. If $\lim_{\rho\rightarrow t,\text{ }0\leq\rho\leq1,\text{ }\xi
\rho/u-\rho+1>0}$ $G_{u}(t)$ exists, we extend $G_{u}(t)$ by the value.
Suppose $\widetilde{G}<+\infty$, then, the equation $G_{u}(t_{u})=\sup_{0\leq
t\leq1,\text{ }\xi t/u-t+1>0}G_{u}(t)$ implies that $t_{u}$ is the optimal proportion.

\bigskip

\noindent\textbf{Theorem 5.1}. $t_{u}=\overline{t}_{u}$ $(u>0)$, \textit{where
}$\overline{t}_{u}$\textit{ is a continuous function defined by }(7)\textit{,
}(8)\textit{ and }(9)\textit{.}

\bigskip

\noindent\textbf{Proof}. (1) If $u\geq E$, $t_{u}=0$.

It should be noted that Lemmas 3.2, 3.3 and equality (10) are valid even if
$u\geq E.$ Using $\partial w_{u}(t)/\partial t<0$ and $w_{u}(0^{+}%
)=E/u-1\leq0$, we have $w_{u}(t)<0$ for each $0<t<\min(1$, $u/(u-\xi))$.
Therefore, we have%
\[
\int_{\rho}^{t}w_{u}(t)dt=\int_{I}\log\frac{a(x)t-ut+u}{a(x)\rho-u\rho
+u}d(F(x))<0
\]
for each $0<\rho<t<\min(1$, $u/(u-\xi))$. This implies that $t_{u}=0.$ The
uniqueness of $t_{u}$ can be easily verified.

\bigskip

(2) If $\xi\leq0$ and $0<u\leq\xi+1/H_{\xi}$, $t_{u}=u/(u-\xi)$.

As shown in (1), it is sufficient to show that $\int_{\rho}^{t}w_{u}(t)dt<0$
for each $0<t$ $<\rho<u/(u-\xi)$. It is clear that $0<H_{\xi}<+\infty.$ Using
$w_{u}(u/(u-\xi)^{-})$ $=(1-\xi/u)H_{\xi}(\xi+1/H_{\xi}-u)\geq0$ (Lemma 3.4)
and $\partial w_{u}(t)/\partial t<0$, we have $w_{u}(t)>0$ for each
$0<t<u/(u-\xi).$ Thus, we obtain the conclusion.

\bigskip

(3) If $\xi>0$ and $\xi<u\leq\xi+1/H_{\xi}$, $t_{u}=1$.

Since $u/(u-\xi)>1,$ we can show that $\int_{\rho}^{t}w_{u}(t)dt<0$ for each
$0<t<\rho<1$ as shown in (2).

\bigskip

(4) If $\xi>0$ and $\xi+1/H_{\xi}<u\leq1/H$, $t_{u}=1$.

From Lemmas 3.6 and 3.8, we have $\widetilde{t}_{u}\geq1$ for each $u\in
(\xi+1/H_{\xi},$ $1/H]$. Therefore, from $w_{u}(\widetilde{t}_{u})=0$ and
Lemma 3.2, we have $w_{u}(t)>0$ for each $0<t<1$. Thus, we obtain the
conclusion as shown in (3).

\bigskip

(5) If $\xi>0$ and $0<u\leq\xi$, $t_{u}=1$.

From $u\leq\xi\leq a(x)$, we have $1\leq a(x)t/u-t+1\leq a(x)\rho/u-\rho+1$
for each $0<t<\rho<1$. Therefore, $\int_{I}\log(\left(  a(x)t/u-t+1\right)
/(a(x)\rho/u-\rho+1))d(F(x))\leq0,$ which implies $t_{u}=1$.

\bigskip

(6) If $\xi>0$ and $1/H<u<E$, $t_{u}=\widetilde{t}_{u}$.

It should be noted that $u/(u-\xi)>1$. It is sufficient to show that
$\int_{\widetilde{t}_{u}}^{t}w_{u}(t)dt$ $<0$ for each $0<t<1$. From Lemmas
3.6, 3.7 and 3.8, we have $0<\widetilde{t}_{u}<1$. Moreover, from
$w_{u}(\widetilde{t}_{u})=0$, we have $w_{u}(t)>0$ for each $0<t<\widetilde
{t}_{u}$ and $w_{u}(t)<0$ for each $\widetilde{t}_{u}<t<1$. Therefore, we
obtain the conclusion.

\bigskip

(7) If $\xi\leq0$ and $\max(0,$ $\xi+1/H_{\xi})<u<E$, $t_{u}=\widetilde{t}%
_{u}$.

It should be noted that $u/(u-\xi)\leq1$. It is sufficient to show that
$\int_{\widetilde{t}_{u}}^{t}w_{u}(t)dt$ $<0$ for each $0<t<u/(u-\xi)$. From
Lemmas 3.2, 3.3 and 3.4, we have $0<\widetilde{t}_{u}$ $<u/(u-\xi)$. Thus, we
obtain the conclusion as shown in (6).\hfill$\square$

\bigskip

Hereafter, we assume that $\widetilde{G}<+\infty.$ Thus, it is easy to verify
the following corollaries.

\bigskip

\noindent\textbf{Corollary 5.1.} \textit{Suppose }$\xi\geq0$\textit{ and
}$1/H<u<E$\textit{, or }$\xi<0$\textit{ and }$\max(0,$\textit{ }$\xi+1/H_{\xi
})<u<E.$\textit{ Then, the optimal proportion of investment }$t_{u}$\textit{
is uniquely determined by }$\int_{I}(a(x)-u)/(a(x)t_{u}-ut_{u}+u)d(F(x))=0$%
\textit{, and the maximized limit expectation of the growth rate is }%
$\exp\left(  \int_{I}\log\left(  a(x)t_{u}/u-t_{u}+1\right)  d(F(x))\right)  $.

\bigskip

\noindent\textbf{Corollary 5.2.} \textit{Suppose }$\xi<0$\textit{ and
}$0<u\leq\xi+1/H_{\xi}$\textit{. Then, the optimal proportion of investment is
}$u/(u-\xi),$\textit{ and the maximized limit expectation of the growth rate
is }$\exp(\int_{I}\log(a(x)-\xi)d(F(x)))/(u-\xi).$

\bigskip

\noindent\textbf{Corollary 5.3.} \textit{Suppose }$\xi\geq0$\textit{ and
}$0<u\leq1/H$\textit{. Then, the optimal proportion of investment is }%
$1,$\textit{ and the maximized limit expectation of the growth rate is }%
$\exp(\int_{I}\log a(x)$ $d(F(x)))/u.$

\bigskip

\noindent\textbf{Corollary 5.4.} \textit{Suppose }$u\geq E.$\textit{ Then, the
optimal proportion of investment is }$0,$\textit{ and the maximized limit
expectation of the growth rate is} $1$.

\bigskip

\bigskip

\noindent{\large \textbf{6. Game pricing}}

\bigskip

In order to determine the price $u$ of the game, we require the riskless
(simple or continuously compounded) interest rate $r>0$ for a period. If
$\xi\geq0$ and $\widetilde{G}$ $<+\infty$, the solution of the equation
$G_{u}(t_{u})=r+1$ (if $r$ is simple) or $G_{u}(t_{u})$ $=e^{r\text{ }}$(if
$r$ is continuously compounded)\ is uniquely determined (Theorems 4.1 and
5.1). In particular, if $u\in(1/H$, $E)$, $t_{u}$ is uniquely determined by
the equation $\int_{I}(a(x)-u)/(a(x)t_{u}-ut_{u}+u)d(F(x))=0$ (Corollary 5.1).
If $\xi<0$, $\widetilde{G}<+\infty$ and $\exp\left(  \int_{I}\log\left(
a(x)\eta+1\right)  d(F(x))\right)  >r+1$ (or $e^{r}$), the solution of the
equation $G_{u}(t_{u})=r+1$ (or $e^{r}$) is uniquely determined (Theorems 4.1
and 5.1).

In this section we assume that $r=0.04$. It is easy to verify that the
following examples are effective games (Lemma 5.3).

\bigskip

\noindent\textbf{Example 6.1. }Suppose that the profit and distribution
functions are given by $a(x)=x$ and\ $F(x)=x\in I=[0$, $1]$ respectively, then
$\xi=0$, $E=1/2,$ and $H=+\infty$. Set $y=t_{u}/u$ $(0<u<1/2)$, then the
equation $w_{u}(t_{u})=0$ can be reduced to$\int_{I}1/(xy$ $-t_{u}+1)dx=1.$
This integral equation has the solution $t_{u}=(e^{y}-y-1)/(e^{y}-1)$.
Therefore, we obtain
\[
G_{u}(t_{u})=\frac{y}{e^{y}-1}\exp\left(  y-1+\frac{y}{e^{y}-1}\right)  ,
\]
which strictly increases from $1$ to $+\infty$ with respect to $y\in(0$,
$+\infty)$. The price $u$ should be the solution of the equation $G_{u}%
(t_{u})=1.04$ (if $r=0.04$ is simple). Thus, the price is $u\fallingdotseq
0.4195$, where $t_{u}\fallingdotseq0.4118$ $(y\fallingdotseq0.9818)$.

It should be noted that the equation $E/u=1.04$ implies the (higher) price
$u\fallingdotseq0.4808$ ($>0.4195$), where $t_{u}\fallingdotseq0.1109$ and
$G_{u}(t_{u})$ $\fallingdotseq1.0022$ ($<1.04$).

\bigskip

\noindent\textbf{Example 6.2. }Suppose that the profit $a$ or $b$ $(a>1>b)$
occurs with probability $p$ or $q=1-p,$ respectively. Further assume that
$1/H<u=1<E$ (if $b>0$) or $u=1<E$\ (if $b\leq0$). Then, from $(a-u)p/(at_{u}%
-ut_{u}+u)+(b-u)q/(bt_{u}-ut_{u}+u)=0$, we obtain
\begin{align}
t_{1}  &  =\frac{q}{1-a}+\frac{p}{1-b}\text{ \ \ (Kelly (1956)), }\tag{17}\\
G_{1}(t_{1})  &  =(a-b)\left(  \frac{q}{a-1}\right)  ^{q}\left(  \frac{p}%
{1-b}\right)  ^{p}.\nonumber
\end{align}

Samuelson (1971) deals with the case in which $a=2.7$, $b=0.3$, and $p=q=0.5$,
where $\xi=0.3$, $E=1.5$, $1/H=0.54$ and $t_{1}=50/119\fallingdotseq0.4202$.
However, Samuelson (1971) may have misinterpreted the criterion to be the
geometric mean $2.7^{0.5}0.3^{0.5}=0.9<1$, instead of $G_{1}(t_{1}%
)=(2.7-0.3)(0.5/1.7)^{0.5}(0.5/0.7)^{0.5}$ $\fallingdotseq1.1000>1$.

When $u\fallingdotseq1.1704$, we have $t_{u}\fallingdotseq0.2898$ and
$G_{u}(t_{u})=1.04$.

It should be noted that the equation $E/u=1.04$ implies the (higher) price
$u\fallingdotseq1.4423$ ($>1.1704$), where $t_{u}\fallingdotseq0.0579,$ and
$G_{u}(t_{u})\fallingdotseq1.0012$ ($<1.04$).

\bigskip

\noindent\textbf{Example 6.3. }In order to obtain the optimal proportions
$(t,$ $s)$ of two independent games $(a(x),$ $F(x))$ and $(b(x),$ $G(x))$ with
the same price $u>0$, we should find the supremum of the function
\begin{equation}
\exp\left(  \int_{J}\int_{I}\log\left(  (a(x)t+b(y)s)/u-t-s+1\right)
d(F(x))d(G(y))\right)  \text{,} \tag{18}%
\end{equation}
where $t\geq0$, $s\geq0$, $s+t\leq1,$ and $(\xi_{a}t+\xi_{b}s)/u-t-s+1>0$. It
should be noted that $\xi_{a}$ (or $\xi_{b}$) is the essential infimum of
$a(x)$ (or $b(x)$).

In Section 5, we introduced the following two games:

\noindent\textbf{Game-1.} The profit $3$ or $0$ occurs with probability $0.5$.
In this case, we have $E$ $=1.5$, $\xi=0,$ $\max_{0\leq t\leq1,\text{ }\xi
t-t+1>0}$ $G_{1}(t)=\sqrt{9/8}\fallingdotseq1.0607,$ and $t_{1}=0.25$. When
$u\fallingdotseq1.0880$, we have $t_{u}\fallingdotseq0.2155$ and $G_{u}%
(t_{u})=1.04$.

\noindent\textbf{Game-2.} The profit $2$ or $0.625$ occurs with probability
$0.5$. In this case, we have $E=1.3125$, $\xi=0.625,$ $\max_{0\leq
t\leq1,\text{ }\xi t-t+1>0}$ $G_{1}(t)=\sqrt{121/96}\fallingdotseq1.1227,$ and
$t_{1}=5/6\fallingdotseq0.8333$. When $u\fallingdotseq1.1237$, we have
$t_{u}\fallingdotseq0.4856$ and $G_{u}(t_{u})=1.04$.

If the same price is $u=1$, the limit expectation of growth rate is given by%
\[
\exp\left(  \left(
\begin{array}
[c]{c}%
\log(3t+2s-t-s+1)+\log(3t+0.625s-t-s+1)\\
+\log(2s-t-s+1)+\log(0.625s-t-s+1)
\end{array}
\right)  /4\right)  \text{.}%
\]
This function attains the maximum value $1.1798$ at $t\fallingdotseq0.2142$
and $s\fallingdotseq0.7809$.

If the same price is $u\fallingdotseq1.8153$, the maximized limit expectation
of growth rate is $r+1=1.04$ at $t\fallingdotseq0.1683$ and $s\fallingdotseq
0.3175$.

\bigskip

\noindent\textbf{Example 6.4. }In the St. Petersburg game (Bernoulli(1738;
English trans. 1954)), suppose that the profit $2^{j}$ occurs with probability
$1/2^{j}$ ($j=1,2,\cdots.$), then $\xi$ $=2$, $E=+\infty,$ and $H=1/3$. This
game is effective, because $\sum_{j=1}^{\infty}(2^{j})^{1/2}/2^{j}<+\infty$.
From Lemma 4.21 we have $G_{1/H}(t_{1/H})=1/3$ $\times\exp(\sum_{j=1}^{\infty
}(\log2^{j})/2^{j})=4/3$. Thus, $G_{u}(t_{u})$ ($u\in(3$, $+\infty)$) strictly
decreases from $4/3$ to $1$. The equation $G_{u}(t_{u})=1.04$ yields the price
$u\fallingdotseq5.1052$. Therefore, if the investors invest $t_{u}%
\fallingdotseq0.1658$\ of their current capital, they can maximize the limit
expectation of growth rate to $1.04.$

It should be noted that the equation $E/u=1.04$ has no solution.

\bigskip

\noindent\textbf{Example 6.5. }The lognormal distributed game is given by%
\begin{equation}
a(x)=Se^{r}e^{x},\text{ \ \ }d(F(x))=\frac{1}{\sqrt{2\pi}\sigma}%
e^{-\frac{(x+\sigma^{2}/2)^{2}}{2\sigma^{2}}}dx, \tag{19}%
\end{equation}
and $I=(-\infty,+\infty)$. In this case, we have $E=Se^{r}$, $H=(e^{-r+\sigma
^{2}})/S,$ and $\exp\left(  \int_{I}\log a(x)d(F(x))\right)  $ $=Se^{r-\sigma
^{2}/2}.$

When $S=100$, $\sigma=0.3,$ and $r=0.04$, we have $\xi=0$, $E\fallingdotseq
104.0811,$ $H$ $\fallingdotseq0.0105125,$ and $1/H\fallingdotseq95.1230.$ From
Lemma 4.16, $G_{u}(t_{u})$\ $(u\in(1/H,$ $E))$ strictly decreases from
$H\exp\left(  \int_{I}\log a(x)d(F(x))\right)  =e^{\sigma^{2}/2}$
$\fallingdotseq1.0460$ to $1$. The equation $G_{u}(t_{u})=e^{0.04}%
\fallingdotseq1.0408$ yields the price $u\fallingdotseq95.6132$. Therefore, if
the investors invest $t_{u}\fallingdotseq0.9450$\ of their current capital,
then they can maximize the limit expectation of growth rate to $e^{0.04}$.

It is clear that the equation $E/u=e^{r}$ yields the (higher) price $u=100$
($>95.6132$). Under this price, if the investors invest $t_{u}\fallingdotseq
0.4433$\ of their current capital, they can maximize the limit expectation of
growth rate to $1.0088$ ($<e^{0.04}$ $\fallingdotseq1.0408$). Because
$\exp(r-\sigma^{2}/2)$ $\fallingdotseq0.9950$ ($<1.0088$), the statement that
the expected growth rate is equal to $r-\sigma^{2}/2$ (Luenberger (1998) 15.5)
is not necessarily true.

\bigskip

\noindent\textbf{Example 6.6. }The European put option is given by%
\begin{equation}
a(x)=\max(K-Se^{rT}e^{x},0),\text{ \ \ }d(F(x))=\frac{1}{\sqrt{2\pi T}\sigma
}e^{-\frac{(x+\sigma^{2}T/2)^{2}}{2\sigma^{2}T}}dx, \tag{20}%
\end{equation}
and $I=(-\infty,+\infty)$. We assume that the stock price $Y=Se^{rT}e^{X}$ is
lognormally distributed with volatility $\sigma\sqrt{T}$, where $S$ is the
current stock price, $r$ is the continuously compounded interest rate, $K$ is
the exercise price of the put option, and $T$ is the exercise period. The
expectation $E$ is given by
\begin{align}
E  &  =\frac{1}{\sqrt{2\pi T}\sigma}\int_{-\infty}^{\log\frac{K}{S}%
-rT}(K-Se^{rT}e^{x})e^{-\frac{(x+\sigma^{2}T/2)^{2}}{2\sigma^{2}T}}%
dx\tag{21}\\
&  =KN\left(  -\frac{\log\frac{S}{K}+(r-\frac{\sigma^{2}}{2})T}{\sigma\sqrt
{T}}\right)  -Se^{rT}N\left(  -\frac{\log\frac{S}{K}+(r+\frac{\sigma^{2}}%
{2})T}{\sigma\sqrt{T}}\right)  ,\nonumber
\end{align}
where $N(x)=\int_{-\infty}^{x}e^{-x^{2}/2}/\sqrt{2\pi}dx$ is the cumulative
standard normal distribution function.

When $S=90$, $K=120$, $T=2$, $\sigma=0.1$, and $r=0.04$, we have $\xi=0$, $E$
$\fallingdotseq22.9848,$ and $H=+\infty$. Therefore, from Theorems 4.1 and
5.1, $G_{u}(t_{u})$\ ($u\in(0,$ $E)$) strictly decreases from $+\infty$ to
$1$. The equations $w_{u}(t_{u})=0$ and $G_{u}(t_{u})=e^{0.08}$ yield the
price $u\fallingdotseq17.8157$. With this price, if investors continue to
invest $t_{u}$ $\fallingdotseq0.5434$\ of their current capital, they can
maximize the limit expectation of growth rate to $e^{0.08}\fallingdotseq
1.0833.$

In general, the equation $E/u=e^{rT}$ yields the price%
\begin{equation}
u=Ke^{-rT}N\left(  -\frac{\log\frac{S}{K}+(r-\frac{\sigma^{2}}{2})T}%
{\sigma\sqrt{T}}\right)  -SN\left(  -\frac{\log\frac{S}{K}+(r+\frac{\sigma
^{2}}{2})T}{\sigma\sqrt{T}}\right)  , \tag{22}%
\end{equation}
which is the Black-Scholes formula for a European put option. Substituting the
above-mentioned values for this formula, we obtain the (higher) price
$u\fallingdotseq21.2176$ ($>17.8157$). With this price, if the investors
continue to invest $t_{u}\fallingdotseq0.2278$ of their current capital, they
can maximize the limit expectation of growth rate to $1.0096$ ($<1.0833$)$.$

\bigskip

\bigskip

\noindent\textbf{References}

\bigskip

\noindent Algoet, P.H., Cover, T., 1988. Asymptotic optimality and asymptotic equipartition

properties of log-optimum investment. The Annals of Probability 16, 876-898.

\noindent Bernoulli, D., 1954. Exposition of a new theory of the measurement of\ risk.

Econometrica 22, 23--36.

\noindent Black, F., Scholes, M., 1973. The pricing of options and corporate liabilities.

Journal of Political Economy 81, 637-654.

\noindent Browne, S., Whitt, W., 1996. Portfolio choice and the Bayesian Kelly criterion.

Advances in Applied Probability 28, 1145-1176.

\noindent Carter, M., van Brunt, B., 2000. The Lebesgue Stieltjes integral. Springer-Verlag,

New York.

\noindent Dutka, J., 1988. On the St. Petersburg paradox. Archive for History
of Exact

\ Sciences 39, 13--39.

\noindent Feller, W., 1957. An introduction to probability theory and \ its
application. John

Wiley and Sons, New York.

\noindent Karatzas, I., Shreve, S.E., 1998. Methods of mathematical finance. Springer-

Verlag, New York.

\noindent Kelly, J.L.Jr., 1956. A new interpretation of information
rate.\ Bell system Technical

Journal 35, 917--926.

\noindent Luenberger, D.G., 1993. A preference foundation for log
mean-variance criteria in

portfolio choice problems. Journal of Economic Dynamics and Control 17,

887-906.

\noindent Luenberger, D.G., 1998. Investment science. Oxford University Press, Oxford.

\noindent Robbins, H., 1961. Recurrent games and the Petersburg\ paradox.
Annals of

Mathematical Statistics 32, 187--194.

\noindent Rotando, L.M., Thorp, E.O., 1992. The Kelly criterion\ and\ the stock\ market.\ 

American Mathematical Monthly 99, 922-931.

\noindent Samuelson, P., 1971. The fallacy of maximizing the geometric mean in long

sequences of investing or gambling. Proceedings of the National Academy of

Sciences of the United States of America 68, 2493--2496.

\noindent Samuelson, P., 1977. St. Petersburg paradoxes: defanged,\ dissected,\ and

historically\ described. Journal of Economic\ Literature 15, 24--55.

\bigskip

\bigskip

\bigskip

\bigskip{\small \noindent\textit{E-mail address}: yukioh@cnc.chukyo-u.ac.jp}

\bigskip

\bigskip

This paper is followed by the articles such as

\noindent Hirashita, Y., 2007. Ratio of price to expectation and complete
Bernstein functions.

Preprint, arXiv:math.OC/0703077.

\noindent Hirashita, Y., 2007, Translation invariance of investment. Preprint, arXiv:math.OC/0703078.

\noindent Hirashita, Y., 2007, Least-Squares prices of games. Preprint, arXiv:math.OC/0703079.

\end{document}